\documentclass[11pt]{article}

\usepackage{arxiv}

\usepackage{amsmath,amsfonts,latexsym,fullpage,graphicx,vmargin, mathrsfs}
\usepackage{pstricks}
\usepackage{comment}
\usepackage[inline]{enumitem}
\usepackage{cleveref}

\usepackage[ruled]{algorithm2e}
\usepackage{algpseudocode}

\RequirePackage[OT1]{fontenc}

\usepackage{amsmath}
\usepackage{amsfonts}
\usepackage{caption}
\usepackage{subcaption}
\usepackage{multirow}
\usepackage{tikz}
\usetikzlibrary{graphs}
\usetikzlibrary{cd}
\usepackage{bm}
\usepackage{url}
\usepackage{placeins}

\RequirePackage{natbib}

\newtheorem{teo}{Theorem}
\newtheorem{prop}{Proposition}
\newtheorem{defi}{Definition}
\newtheorem{rmk}{Remark}
\newtheorem{lem}{Lemma}
\newtheorem{cor}{Corollary}
\newtheorem{assump}{Assumption}
\newtheorem{claim}{Claim}

\DeclareMathOperator*{\LCA}{LCA}

\DeclareMathOperator*{\len}{len}
\DeclareMathOperator*{\lvl}{lvl}

\usepackage{scalerel}

\newcommand{\R}{\mathbb{R}}

\newcommand{\virgolette}[1]{``#1''}

\firstpageno{1}

\begin{document}

\title{A Graph-Matching Formulation of the Interleaving Distance between Merge Trees}

\author{Matteo Pegoraro\thanks{Department of Mathematical Sciences, Aalborg University}}

%
\maketitle

\begin{abstract}
In this work we study the interleaving distance between merge trees from a combinatorial point of view. We use a particular type of matching between trees to obtain a novel formulation of the distance. With such formulation, we tackle the problem of approximating the interleaving distance by solving linear binary optimization problems in a recursive and dynamical fashion, obtaining lower and upper bounds. We implement those algorithms to compare the outputs with another approximation procedure presented by other authors. We believe that further research in this direction could lead to polynomial time algorithms to approximate the distance and novel theoretical developments on the topic.
\end{abstract}


\section{Introduction}

Topological data analysis (TDA)
is a scientific field lying at the crossroads of topology and data analysis: objects like functions and point clouds are typically studied by means of - possibly multidimensional - complexes of homology groups \citep{hatcher} obtained with different pipelines. 
Homology groups - considered with coefficients in a field $\mathbb{K}$ - are vector spaces, whose dimension is determined by different kind of \virgolette{holes} which can be found in a space and thus provide a rich characterization of the shape of the space the analyst is considering.
More precisely, each statistical unit induces a whole family of topological spaces indexed on $\R^n$
which then, via homology, produces a family of vector 
spaces indexed on the same set - more precisely, a \emph{functor} \citep{maclane} $(\mathbb{R}^n, \leq )\rightarrow \text{Vect}_\mathbb{K}$.
Such collections of vector spaces are called \emph{persistence modules} \citep{chazal2008persistencemodules} - with \emph{multidimensional} persistence modules being a special reference to the cases in which $n>1$ - and describe the shape of the data by considering interpretable topological information at different resolutions. 
A topological summary or an invariant of a persistence module is a representation which usually maps persistence modules into a metric space - or even a vector space, so that some kind of analysis can be carried out.
 The most used topological summaries for 1-D persistence (i.e. $n=1$) include persistence diagrams \citep{PD_1}, persistence landscapes \citep{landscapes}, persistence images \citep{pers_img} and persistence silhouettes \citep{silhouettes}. 
When dealing with $0$-dimensional homology and 1-D persistence and with a space $X$ which is path connected, however, another topological summary called \emph{merge tree} \citep{merge_interl} can be employed.
A merge tree is a tree-shaped topological summary which captures the evolution of the path connected components $\pi_0(X_t)$ of a filtration of topological spaces $\{X_t\}_{t\in\R}$, with $X_t\hookrightarrow X_{t'}$ if $t\leq t'$. 


One of the central ideas in TDA is the one of \emph{stability} properties, i.e. the ones related to the continuity of the operator mapping data into topological signature. One definition which has gained a lot of success is the one of \emph{$\varepsilon$-interleavings} between topological representations: when adding some kind of $\varepsilon$-noise to the data, the associated summary \virgolette{moves} by at most $\varepsilon$.
If this happens, then the operator mapping functions to topological representations is a non-expansive Lipschitz operator.
Starting from the \emph{bottleneck} distance for \emph{persistence diagrams} \citep{PH_survey}, this idea has been applied to in 1-D persistence modules \citep{chazal2008persistencemodules}, multidimensional persistence modules \citep{interl_multi}, merge trees \citep{merge_interl}, sheaves \citep{munch2016convergence, curry2023convergence} and more general situations \citep{de2017theory, interl_homotopy}.
In this paper we focus on the problem of studying interleavings for merge trees.

The main reason to pursue such research direction is that merge trees and persistence diagrams are not equivalent as topological summaries, that is, they do not represent the same information with merge trees being able to distinguish between situations which persistence diagrams - and all the other equivalent summaries - cannot \citep{kanari2020trees, curry2021decorated, mergegrams, smith2022families}.
As a consequence, in recent years, a lot of research sparkled on merge trees, mainly driven by the 
need of having a \emph{stable} metric structure to compare such objects. Stable and unstable edit distances between Reeb graphs or merge trees have been proposed by \cite{merge_farlocca, di2016edit, bauer2020reeb,pegoraro2024finitely}, while other works focus on (unstable) Wasserstein Distances \citep{merge_wass}, $l_p$ distances
\citep{cardona2021universal} and  interleaving distances between merge trees
\citep{merge_interl, interl_approx_2018, interl_approx, interl_np_hard,  merge_intrins} or Reeb graphs 
\citep{bauer2015strong, de2016categorified, munch2016convergence, curry2023convergence}.

\subsection*{Related Works}

The problem of computing the interleaving distance between merge trees has already been approached by \cite{interl_approx_2018}
and \cite{interl_approx} in an attempt to approximate the Gromov-Hausdorff (GH) distance when both metric spaces are metric trees. In this situation, in fact, up to carefully choosing the root for the metric trees, the two metrics are equivalent.
Both problems - approximating the GH distance and computing the interleaving distance - have been shown to be NP-hard \citep{interl_approx_2018, interl_np_hard} and thus even obtaining feasible approximation algorithms for small trees is a daunting task.

In \cite{interl_approx_2018} the authors provide an algorithm to approximate the interleaving distance between merge trees via binary search. They build a set which contains the solution of the optimization problem and then obtain a criterion to assess if certain values of the aforementioned set can be excluded from the set of solutions. If the length of the branches of both trees is big enough, one can carry out the binary search over all possible couples of vertices, with one vertex from the first tree and one from the second.  
If the trees instead posses branches which are too small, the decision procedure is coupled with a \emph{trimming} step, to deal with such smaller edges. The algorithm returns an approximation of the interleaving distance, with the approximation factor depending on the ratio between the longest and the shortest edge in the tree and the number of vertices.

\cite{interl_approx}, instead, propose the first algorithm for the exact computation of the interleaving distance. Starting from a novel definition of the interleaving distance, they develop a faster alternative to exclude values from the same set of solutions considered by \cite{interl_approx_2018}.
By filtering a finite subset of points on the metric trees - in general bigger then the set of vertices of the \virgolette{combinatorial} trees - according to a) their heights b) a candidate optimal value $\delta$, in an up-bottom fashion, points of the second tree are matched to subsets of points of the first tree. Such matching means that all points of the first tree inside in the chosen collection can be potentially mapped in the point of the second tree via some function viable for the interleaving distance. Each such matching $(S,w)$ is then used to establish similar couples - \emph{(set, point)} - between the children of the involved points: the set of points of the first tree is to be chosen among the compatible subsets of the children of $S$, and the point on the second tree among the children of $w$. If no such couples are found $\delta$ is discarded. 
The algorithm they propose to compute the interleaving distance has complexity $O(n^2 \log^3(n) 2^{2\tau} \tau^{\tau + 2})$ where $n$ is the sum of the vertices in the two merge trees, and $\tau$ is a parameter - depending on the input trees - which is defined by the authors. The key issue is that $\tau$ can be very big also for small-sized trees: in \cite{interl_approx}, Figure 2, the authors showcase a tree with $8$ leaves and $\tau = 13$ (thus the complexity of the algorithm is at least $O(n^2 \log^3(n) 2^{26} 13^{15})  \sim O(n^2 \log^3(n) 10^{24}$).

Due to its computational complexity the algorithm by \cite{interl_approx} has not been implemented by other authors working with merge trees \citep{curry2021decorated}, which instead have exploited another formulation of the same metric, provided by \cite{merge_intrins}.
This last formulation relies on a definition of the interleaving distance between merge trees which are endowed with a fixed set of labels on their vertices. The usual interleaving distance is then shown to be equivalent to choosing an appropriate set of labels - potentially upon adding some vertices to the tree - for the two given merge trees. This formulation, per se, does not provide computational advantages over the classical one, since evaluating all the possible labelings of a tree would still be unfeasible. However, 
\cite{curry2021decorated} propose a labeling strategy which should provide \emph{good} labels. Despite being computationally accessible even for very big trees, this approach has the downside of not providing methods for assessing the goodness of the approximation.

\subsection*{Contributions}

As already mentioned, in the present paper we take a perspective which differs from all the aforementioned works: instead of obtaining progressively sharper bounds for the distance by locally looking at differences between trees or instead of looking for optimal labelings, we aim at describing jointly a global matching between two merge trees whose \emph{cost} returns the exact interleaving distance between them. Describing matchings between trees is often useful to formulate optimization problems via binary programming algorithms, which have proven to be a very useful tool to deal with the computational burden introduced by unlabeled or unordered trees \citep{TED}. 
Indeed, we are able to produce a formulation of the interleaving distance originating a mixed binary linear programming (MBLP) approach which gives lower and upper bounds for the interleaving distance between two trees. The computational complexity of such bounds is close to the one of the formulation of the classical tree edit distance (for unordered trees) via binary linear programming (BLP) \cite{TED}. See \Cref{prop:complexity} and \Cref{rmk:comp_cost}. 
 We use the upper bound to assess the approximation method proposed by \cite{curry2021decorated}. 
 Lastly, the formulation we obtain in the present work is used by \cite{pegoraro2024finitely} to obtain some inequalities between the edit distance therein defined and the interleaving distance.

As discussed in the end of \Cref{sec:map_couplings}, \Cref{rmk:comp_cost} and in \Cref{sec:discussion_interl}, we also think that the results in the manuscript can lead to a) a polynomial time approximation of the interleaving distance, following the approach of \cite{zhang1996constrained} b) another formulation of the distance, potentially opening up the door to novel theoretical developments.

\subsection*{Outline}

The paper is organized as follows. In \Cref{sec:merge_trees_interl} we introduce merge trees first in a combinatorial fashion and then as \virgolette{continuous} posets and metric spaces, recalling the definition of interleaving distance. In \Cref{sec:couplings} we introduce a partial matching between trees called \emph{coupling}. Relationships between these matchings and maps between trees are then presented in \Cref{sec:coup_maps_costs}.
\Cref{sec:couplings_maps} and \Cref{sec:map_couplings} prove the equivalence between the usual definition of interleaving distance and the problem of finding optimal couplings.
In the last part of the manuscript we face the problem of approximating an optimal matching between trees: in \Cref{sec:prop} we prove two properties of the interleaving distance, which allow us to
write down the dynamical programming algorithm of \Cref{sec:LIP}, while
\Cref{sec:errors} points out few facts about the error propagation in the previously presented algorithm. In \Cref{sec:simulation} we test our approximation with the one proposed by other authors. \Cref{sec:discussion_interl} concludes the main body of the work with a brief discussion.
\Cref{sec:proofs} contains some of the proofs of the results in the manuscript, while the other ones are reported right after the statement they relate to, as they can help the reader in following the discussion.

\section{Merge Trees and Interleaving Distance}\label{sec:merge_trees_interl}

We start by introducing merge trees as combinatorial objects along with their \virgolette{continuous} counterpart.

\subsection{Merge Trees as finite graphs}

In accordance with other authors, we call \emph{merge tree}  a rooted tree with a suitable height function defined on its vertices. 
We now state the formal definition, introducing some pieces of notation which we will use to work with those objects throughout the paper.

\begin{defi}
A tree structure $T$ is given by a set of vertices $V_T$ and a set of edges $E_T\subset V_T\times V_T$ which form a connected rooted acyclic graph.  We indicate the root of the tree with  $r_T$. We say that \(T\) is finite if \(V_T\) is finite. The order of a vertex $v \in V_T$ is the number of edges which have that vertex as one of the extremes, and is called $ord_T(v)$. 
Any vertex with an edge connecting it to the root is its child and the root is its father: this is the first step of a recursion which defines the father and children relationship for all vertices in \(V_T.\)
The vertices with no children are called leaves  or taxa and are collected in the set $L_T$. The relation $child < father $ generates a partial order on $V_T$. The edges in $E_T$ are represented by ordered couples $(a,b)$ with $a<b$.
A subtree of a vertex $v$, called $sub_T(v)$, is the tree structure whose set of vertices is $\{x \in V_T\mid x\leq v\}$. 
\end{defi}

Note that, identifying an edge $(v,v')$ with its lower vertex $v$, gives a bijection between $V_T-\{r_T\}$ and $E_T$, that is $E_T\simeq V_T$ as sets. 
Given this bijection, we may use $E_T$ to indicate the vertices $v\in V_T-\{r_T\}$,  to simplify the notation.

Now we give the notion of \emph{least common ancestor} between vertices.

\begin{defi}
Given a set of vertices $A=\{a_1,\ldots,a_n\}\subset V_T$, we define $\LCA(a_1,\ldots,a_n)=\min\bigcap^n_{i=1} \{v\in V_T \mid v\geq a_i\}$.
\end{defi}

Now, to obtain a merge tree we add an increasing height function to a tree structure.

\begin{defi}\label{defi:MT}
A merge tree $(T,f)$ is a finite tree structure $T$ such that the root is of order $>1$, coupled with a monotone increasing function \(f:V_T\rightarrow \mathbb{R}\).
\end{defi}

\begin{rmk}
\Cref{defi:MT} is slightly different from the definition of merge trees found in \cite{merge_intrins}, \cite{pegoraro2024finitelyfunc} and other works. In particular, in the present work, we do not need to have a root at infinity and thus we remove it to avoid unpleasant technicalities. 
Similarly, the function coupled with the tree structure in literature is usually referred to as $h_T$ being an \virgolette{height} function. To avoid ovearloading the notation, since we need to introduce many subscript and superscripts, we call these functions with more usual functional notations like $f$ or $g$. 
\end{rmk}

\begin{assump}
To avoid formal complications we make the following standard genericity assumption for any merge tree $(T,f)$:
\begin{itemize}
\item[(G)] $f:V_T\rightarrow \mathbb{R}$ is injective.
\end{itemize}
\end{assump}

\begin{example}
Given a real valued Morse \citep{milnor1963morse} function $f:X\rightarrow \mathbb{R}$ defined on a path connected compact space $X$, the \emph{sublevel set} filtration is given by $X_{t}=f^{-1}((-\infty,t])$, together with the maps $X_{t<t'}=i:f^{-1}((-\infty,t])\hookrightarrow f^{-1}((-\infty,t'])$.
We can describe the family of sets and maps $\{\pi_0(X_t)\}_{t\in \R}$ via a merge tree \citep{pegoraro2024finitelyfunc}.
Similarly we can consider a finite set $C\subset \mathbb{R}^n$ and take its  \emph{C\'ech} filtration: $X_{t}=\bigcup_{c\in C}B_{t}(c) $. With $B_{t}(c)=\{x\in\mathbb{R}^n \mid \parallel c-x \parallel < t\}$. As before: $X_{t<t'}=i:\bigcup_{c\in C}B_{t}(c)\hookrightarrow \bigcup_{c\in C}B_{t'}(c)$ and
$\{\pi_0(X_t)\}_{t\in \R}$ can be represented via a merge tree \citep{pegoraro2024finitelyfunc}.
\end{example}

Before proceeding we need one last graph-related definition which we use to denote some particular kinds of paths on the tree structure of a merge tree.

\begin{defi}
Given a merge tree $T$, a sequence of edges is an ordered sequence of adjacent edges $\{e_1,\ldots,e_n\}$. Which means that we have $e_1<\ldots < e_n$, according to the order induced by the bijection $E_T\leftrightarrow V_T-\{r_T\}$ and that $e_i$ and $e_{i+1}$ share a vertex. We will use the notation $[v,v']$ to indicate a sequence of edges which starts in the vertex $v$ and ends before $v'$, with $v<v'$. Thus, $v$ is included in the sequence while $v'$ is the first excluded vertex, coherently with $E_T\simeq V_T-\{r_T\}$.
\end{defi}

In the left column of \Cref{fig:metric_vs_merge}
the reader can find an example of a merge tree, which can be used to get familiar with the definitions just given. For instance one can check that, using the notation of the figure, $R=\LCA(A,B,D)$ and $[A,C,R]$ is an ordered sequence of edges, while $[C,R,D]$ is not.

\subsection{Metric Merge Trees}

\begin{figure}
    \begin{subfigure}[c]{0.97\textwidth}
    	\centering
    	\includegraphics[width = \textwidth]{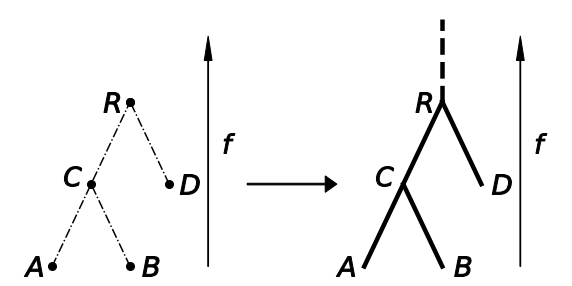}
	\end{subfigure}
\caption{A merge tree (left) with its associated metric merge tree (right).}
\label{fig:metric_vs_merge}
\end{figure}

Now we consider the continuous version of a merge tree, intuitively obtained by considering all the points in the edges of a merge tree as points of the tree itself. 
For a visual intuition of the following ideas, the reader may refer to 
\Cref{fig:metric_vs_merge}. 

\begin{defi}
Given a merge tree $T$, we obtain the associated metric merge tree as follows:
\[
\mathbf{T}=[f(r_T),+\infty )\coprod_{(x,x')\in E_T} [f(x),f(x')]/\sim 
\]
where, for every $v\in V_T$, $f(v)\sim f(v')$ if $v=v'$. 
We refer to the points in $\mathbf{T}$ with the same notation we use for vertices in $v_T$: 
$x\in\mathbf{T}$; with $f(x)$ we indicate the height values of $x$.

For every point $x\in V_T$, we can identify $x$ with the point $f(x)\in[f(x),f(x')]$, with $(x,x')\in E_T$, or, equivalently,  with  
$f(x)\in[f(x'),f(x)]$, if $(x',x)\in E_T$. This induces a well defined map $V_T\hookrightarrow \mathbf{T}$.
Thus, given $v\in V_T$, with an abuse of notation, we can consider $v\in\mathbf{T}$. Note that the height function $f:\mathbf{T}\rightarrow \R$ extends the function defined on the vertices of the merge tree. Every point in 
$\mathbf{T}-V_T$ belongs to one and only one interval of the form $[f(x),f(x')]$. Thus we can induce a partial order relationship on $\mathbf{T}$ by ordering points first with the partial order in $E_T$ and then with the increasing internal order of $[f(x),f(x')]$. Thus we can explicitly write down the shortest path metric in $\mathbf{T}$: $ d(x,x')=f(\LCA(x,x'))-f(x)+f(\LCA(x,x'))-f(x')$. 
\end{defi}

\begin{rmk}[Notation and Previous Works]
If we start from a morse function $f:D\rightarrow \R$ defined on a compact and path connected topological space $D$, what we call \virgolette{metric merge tree} is essentially
the \virgolette{merge tree} of $f$, according to the topological definition given in \cite{merge_interl}, also called \emph{classical merge tree} in \cite{curry2021decorated}. 
Also \emph{display posets} of persistent sets, as defined in \cite{curry2021decorated}, are - up to some additional hypotheses - equivalent to metric merge trees. But 1) as most of the previous literature on such topic, we are interested at working on merge trees in a combinatorial fashion, so we build merge trees as graphs and then bridge to this continuous construction 2) to use display posets we would need to give other technical definitions which we do not need in this work, thus we avoid working with such objects. To sum up, all these procedures/definitions are different ways to build an $\R$-space $(X,f_X:X\rightarrow \R)$ - with $f_X$ continuous function - \citep{de2016categorified} with the underlying topological space $X$ being the CW-complex associated to a tree.
In particular the display poset highlights the poset structure of such topological space while 
the approach in \cite{merge_interl} describes a geometric construction which gives the display poset a natural topology. Our metric merge tree in some sense condenses all these properties in one object. 

Lastly, starting from any of the aforementioned constructions leading to the poset and $\R$-space 
$(\mathbf{T},f:\mathbf{T}\rightarrow \R)$ the discrete object we call merge tree can easily be obtained using the (unique) transitive reduction \citep{transitive_reduction} of the directed acyclic graph associated to (a finite representation of) the display poset of $\mathbf{T}$, with the induced height function - see also \citep{pegoraro2023edit}.
\end{rmk}

For each metric tree $\mathbf{T}$, we have 
a family of continuous maps, $s^k_T:\mathbf{T}\rightarrow \mathbf{T}$, with $k\geq0$, called \emph{structural maps}, defined as follows: $s^k_T(x)=x'$ with $x'$ being the only point in $\mathbf{T}$ such that $x'\geq x$ and $f(x')=f(x)+k$.

\subsection{Interleaving Distance between Merge Trees}

Now we recall the main facts about the interleaving distance between merge trees, relying on \citep{merge_interl}.

\begin{figure}
    \begin{subfigure}[c]{0.97\textwidth}
    	\centering
    	\includegraphics[width = \textwidth]{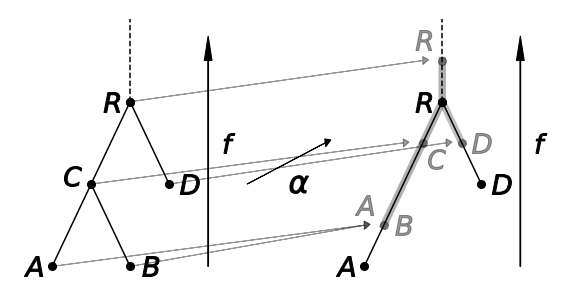}
		\caption{A graphical representation of a continuous function $\alpha:\mathbf{T}\rightarrow \mathbf{G}$ between metric merge trees satisfying condition \emph{(I1)} in \Cref{def:interleaving}.}
		\label{fig:alpha}
	\end{subfigure}

    \begin{subfigure}[c]{0.97\textwidth}
    	\centering
    	\includegraphics[width = \textwidth]{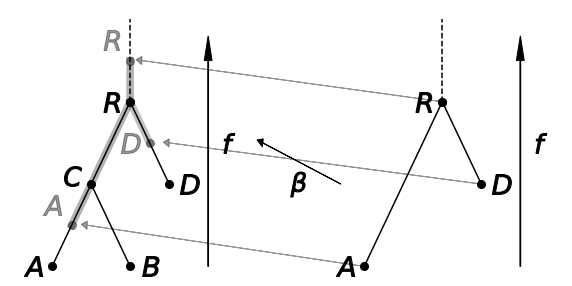}
		\caption{A graphical representation of a continuous function $\beta:\mathbf{G}\rightarrow \mathbf{T}$ between metric merge trees satisfying condition \emph{(I1)} in \Cref{def:interleaving}.}
		\label{fig:beta}
	\end{subfigure}	
\caption{An example of maps $\alpha$ and $\beta$ giving the $\varepsilon$-interleaving of two metric merge trees.}
\label{fig:interleaving}
\end{figure}

\begin{defi}[\cite{merge_interl}]\label{def:interleaving}
Two continuous maps $\alpha:\mathbf{T}\rightarrow \mathbf{G}$ and $\beta:\mathbf{G}\rightarrow \mathbf{T}$ between two metric merge trees $\mathbf{T}$ and $\mathbf{G}$, 
are $\varepsilon$-compatible, with $\varepsilon\geq 0$, if:
\begin{itemize}
\item[(I1)] $g(\alpha(x))=f(x)+\varepsilon$ for all $x\in\mathbf{T}$ and $f(\beta(y))=g(y)+\varepsilon$ for all $y\in\mathbf{G}$;
\item[(I2)] $\alpha\beta = s_G^{2\varepsilon}$ and $\beta\alpha = s_T^{2\varepsilon}$.
\end{itemize} 
The interleaving distance between $\mathbf{T}$ and $\mathbf{G}$ is then: $ d_I(\mathbf{T},\mathbf{G})= \inf \{\varepsilon \mid \text{ there 
$\varepsilon$-compatible maps}\}$.
\end{defi}

\begin{figure}
    \begin{subfigure}[c]{0.97\textwidth}
    	\centering
    	\includegraphics[width = \textwidth]{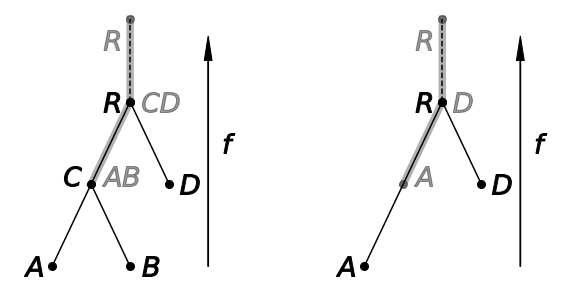}
	\end{subfigure}
\caption{A representation of the images of the maps $\beta\alpha$ (left) and $\alpha\beta$ (right), with $\alpha$ and $\beta$ being the maps in \Cref{fig:interleaving}.}
\label{fig:alpha_beta}
\end{figure}

For an example of $\alpha$ and $\beta$ continuous map satisfying \emph{(I1)} see \Cref{fig:interleaving}. Those maps satisfy also \emph{(I2)} as shown by \Cref{fig:alpha_beta}.

The work of \cite{interl_approx} shows that the existence of $\alpha$ and $\beta$ is in fact equivalent to the existence of a single map $\alpha$
with some additional properties which are stated in the next definition.

\begin{defi}[\cite{interl_approx}]
Given two metric merge trees $\mathbf{T}$ and $\mathbf{G}$, an $\varepsilon$-good map $\alpha:\mathbf{T}\rightarrow \mathbf{G}$ is a continuous map such that:
\begin{itemize}
\item[(P1)] $g(\alpha(x))=f(x)+\varepsilon$ for all $x\in\mathbf{T}$
\item[(P2)] if $\alpha(x)<\alpha(x')$ then $s_T^{2\varepsilon}(x)<s_T^{2\varepsilon}(x')$
\item[(P3)] if $y\in\mathbf{G}-\alpha(\mathbf{T})$, then, given $w = \min \{y'>y \mid y'\in\alpha(\mathbf{T})\}$, we have $g(w)-g(y)\leq 2\varepsilon$. 
\end{itemize}
\end{defi}

As anticipated, \cite{interl_approx} prove that two merge trees are $\varepsilon$-interleaved if and only if there is a $\varepsilon$-good map between them.

\section{Couplings}
\label{sec:couplings}

Now we start our combinatorial investigation. 
Given two merge trees $T$ and $G$, we would like to match their vertices and compute the interleaving distance relying only on such matches.
To match the two graphs we will use a set $C\subset V_T \times V_{G}$ which will tell us which vertex is coupled with which.
Clearly this set $C$ must satisfy some constraints which we now introduce.
 
Given finite sets $A,B$, we indicate with $\#A$ the cardinality of one set and with $A-B:= \{a\in A \mid a\notin B\}$. Consider now $C\subset V_T \times V_{G}$ and the projection $\pi_T:V_T \times V_{G} \rightarrow V_T$, we define the multivalued map $\Lambda_C^T:V_T\rightarrow V_T$ as follows: 
\[ 
\Lambda_C^T(v) =
\begin{cases} 
      \max_{v'<v} \pi_T(C) & \text{if }\#\{ v'\in V_T \mid v'<v \text{ and } v'\in \pi_T(C)\}>0 \\
      \emptyset & \text{otherwise.} 
   \end{cases}
\]
The term \virgolette{multivalued} means that $\Lambda_C^T$ is a function $\Lambda_C^T:V_T\rightarrow \mathcal{P}(V_T)$ with $\mathcal{P}(V_T)$ being the power set of $V_T$.

Since $V_T$ and $V_G$ are posets, we can introduce a partial order relationship on any $C\subset V_T\times V_G$, having $(a,b)<(c,d)$ if and only if $a<c$ and $b<d$.

\begin{defi}
A coupling between two merge trees $(T,f)$ and $(G,g)$ is a set $C\subset V_T \times V_{G}$ such that:
\begin{itemize}
\item[(C1)] $\#\max C=1$ or, equivalently, $\#\max \pi_T(C)=\#\max \pi_G(C)=1$ 
\item[(C2)] the projection $\pi_T:C \rightarrow V_T$ is injective (the same for $G$)
\item[(C3)] given $(a,b)$ and $(c,d)$ in $C$, then $a<c$ if and only if $b<d$
\item[(C4)] $a\in \pi_T(C)$ implies $\#\Lambda_C^T(a)\neq 1$ (the same for $G$). 
\end{itemize}

The set of couplings between $T$ and $G$ is called $\mathcal{C}(T,G)$.
\end{defi}

\begin{figure}
    \begin{subfigure}[c]{0.97\textwidth}
    	\centering
    	\includegraphics[width = \textwidth]{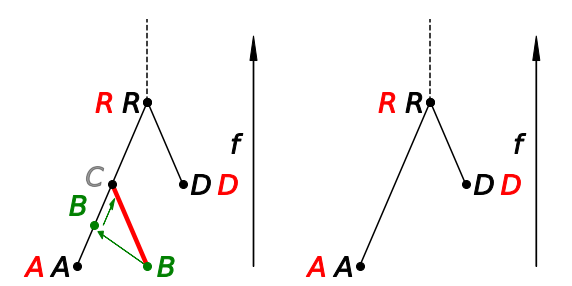}
		\caption{Call $T$ the merge tree on the left and $G$ the one on the right. The couples of adjacent letters - red and black - indicate a coupling. The vertex $C$ in grey represent the only vertex with $\#\Lambda(C)=1$ and belongs to $U^T$; instead the leaf $B$ belongs to $D^T$ and clearly $\#\Lambda(B)=0$. The green arrows suggest the path of the two-step deletion of $B$ - note that $\varphi(B)=C$ and so $\eta(B)=A$.}
		\label{fig:coup_1}
	\end{subfigure}

    \begin{subfigure}[c]{0.97\textwidth}
    	\centering
    	\includegraphics[width = \textwidth]{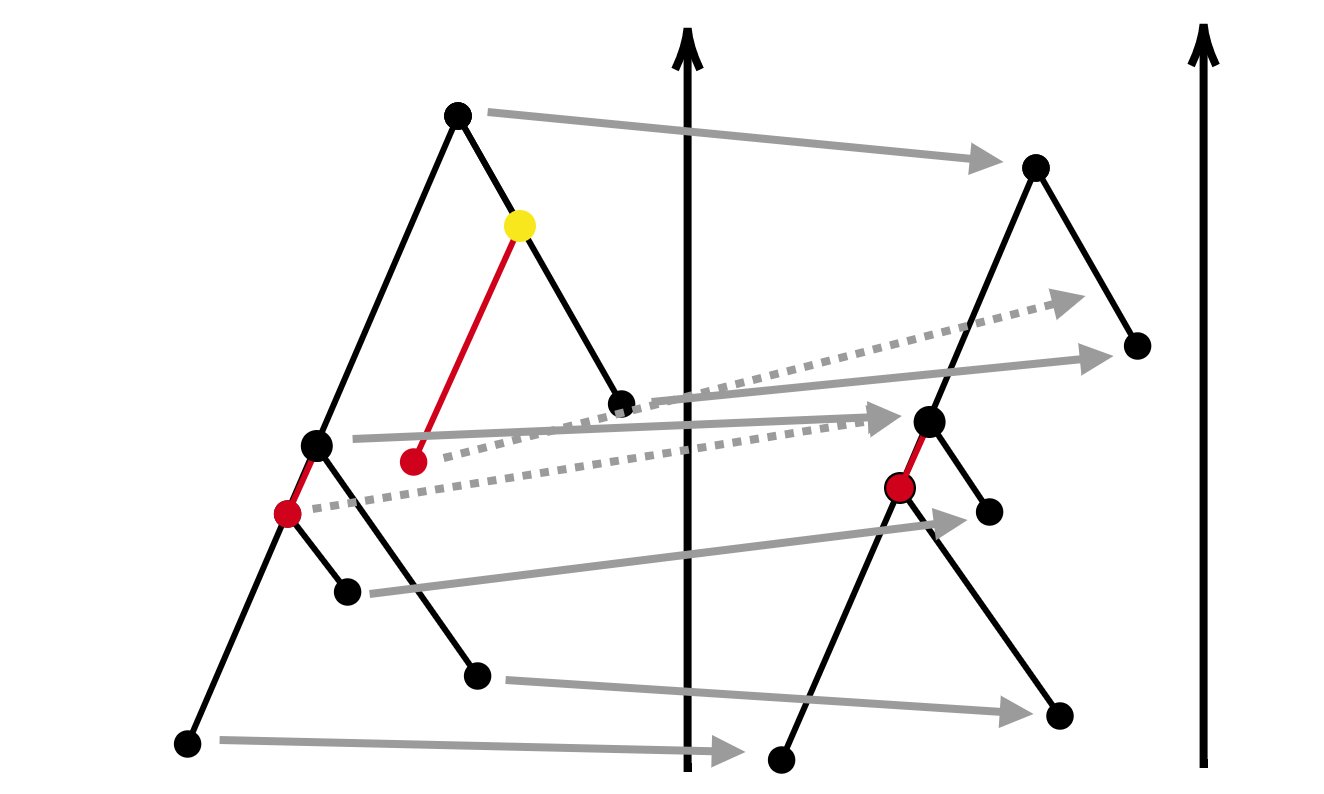}
		\caption{		
		The grey arrows between the two merge trees indicate a coupling $C$. The dashed gray arrow, instead, indicate where the lower vertex of the deleted edges - the red ones - are sent by the associated function $\alpha_C$ - introduced in \Cref{sec:coup_maps_costs}. The lower deleted edge - of the leftmost tree - is an internal edge and its lower vertex $v$ is deleted with $\#\Lambda(v)>1$ and thus it is sent to $\chi(v)$. The higher deleted edge of the leftmost tree is instead deleted with a two step deletion as in \Cref{fig:coup_1}. The yellow vertex is unused.}
		\label{fig:coup_2}
	\end{subfigure}	
\caption{Two examples of couplings, displaying all possible cases of unused vertices and deletions.}
\label{fig:couplings}
\end{figure}

We invite the reader to follow the remaining of the section looking at \Cref{fig:couplings} which can help understanding the comments we present on properties (C1)-(C4) and visualizing the pieces of notation we need to establish.

We collect the set of points $v\in E_T$ such that  $\# \Lambda_C^T(v)=1$ in a set we call $U_C^T$. Instead, the points such that $v\notin\pi_T(C)$ and  $\# \Lambda_C^T(v)\neq 1$ are collected in a set $D_C^T$. 
Note that the sets $\pi_T(C),U_C^T$ and $D_C^T$ are a partition of $V_T$. 
A couple $(v,w)\in C$ means that we match  a vertex $v\in V_T$ with a vertex $w\in V_{G}$. It is clear from the definition that there can be vertices left out from this matching. We regard this vertices as unnecessary for the coupling $C$: when we will induce maps between metric merge trees starting from couplings, the position of these vertices inside the metric merge tree $\mathbf{G}$ will be completely induced by other vertices in $V_T$. Among the vertices not appearing in the couples of $C$ we distinguish between two situations: $v\in D_C^T$ can be informally thought of $(v,father(v))$ needing to be removed the image of $\beta(\alpha(\mathbf{T}))$, with $\alpha$ and $\beta$ as in \Cref{def:interleaving}. For this reason we say that the vertices in $D_C^T$ and $D_C^{G}$ are \emph{deleted} by the coupling $C$. Instead, the vertices $v\in U_C^T$ are vertices which are \emph{unused}, ignored by the coupling, and will be of no importance in the computation of the distance.

In this context, property (C1) is asking that the vertices of $T$ and of $G$ coupled by $C$ form a rooted tree, while (C2) is asking that each vertex is paired with at most one other vertex of the other merge tree. Note the maps $\alpha:\mathbf{T}\rightarrow \mathbf{G}$ in \Cref{def:interleaving} are not forced to be injective, but, as we will prove in later sections, for our purposes it is enough to couple points just one time. Condition (C3) is asking that the coupling respects the tree structures of $T$ and $G$; in particular this implies that $(\max \pi_T(C),\max \pi_G(C))\in C$.
Lastly, due to condition (C4), we avoid coupling vertices which have only one element below them which is coupled. 

Given a coupling, we want to associate to a vertex a cost which indicates how much the coupling moves that vertex.
In order to do so, we define the following functions:

\begin{itemize}
\item define $\varphi^C_T:V_T\rightarrow V_T$ so that 
 $\varphi^C_T(x)=\min \{v\in V_T \mid v> x \text{ and }\#\Lambda(v)\neq 0\}$.
Note that since the set 
$\{v\in V_T\mid v> x\}$ is totally ordered, $\varphi^C_T(x)$ is well defined;
\item similarly, define $\delta^C_T:V_T\rightarrow V_T$, defined as $\delta^C_T(x) =
 \min\{v \in V_T \mid v\geq x \text{ and } v \in \pi_T(C) \}$;
\item define $\chi^C_T:V_T\rightarrow V_G$, defined as $\chi^C_T(x) =
\LCA(\{\pi_G((v,w))\mid v \in \Lambda_T(x)\})$;
\item set $\gamma^C_T:V_T-D_C^T\rightarrow V_G$ to be:
\[ 
\gamma^C_T(x) =
\begin{cases} 
      \arg\min\{g(w) \mid (v,w) \in C, v < x \} & \text{if }\#\{g(w) \mid (v,w) \in C, v < x \}>0 \\
      \emptyset & \text{otherwise.} 
   \end{cases}
\]
Note that if $\#\{g(w) \mid (v,w) \in C, v < x \}>0$, by \emph{(G)}, $\gamma^C_T(x)$  is uniquely defined. Moreover, $\gamma^C_T(\varphi^C_T(x))$ is well defined for any $v\in V_T$. 
\end{itemize}

\begin{rmk}
When clear from the context, to lighten the notation, we might omit the subscripts and superscripts. For instance if we fix $C\in\mathcal{C}(T,G)$ and $x\in V_T$, we refer to $\gamma^C_T(\varphi^C_T(x))$ as $\gamma(\varphi(x))$. And similarly we use 
$x\in U$, $x\in D$ respectively for $x\in U_C^T$ and $x\in D_C^T$.
\end{rmk}

Given $x\in V_T$, in what follows we will encounter the following
objects (see \Cref{fig:couplings}): 
\begin{itemize}
\item  $\eta(x):=\gamma(\varphi(x))$ is the vertex in $V_G$ obtained in this way: starting from $x$ we go towards $r_T$ until we meet the first point of $V_T$ which has at least one vertex in its subtree which is not deleted; we consider the subtree of such vertex and take all the vertices of $V_G$ which are coupled with elements in this subtree (such set is not empty); lastly we take the lowest vertex among this set;
\item given $x$ such that $x\in D$ and $\#\Lambda(x)>1$, we often consider $\delta(x)$; note that for every $v$ such that $x\leq v < \delta(x)$ then $v\in D$ and $\#\Lambda(v)>1$. Thus $[x,\delta(x)]$ contains only deleted vertices, up to $\delta(x)$;
\item similarly, for $x$ such that $x\in D$ and $\#\Lambda(x)>1$, we take $\chi(x)$: that is the the lowest point in $G$ where $x$ can be sent compatibly with $C$ and the tree structures of $T$ and $G$. Thus, if $(x,y)\in C$, $y\geq \chi(x)$.	  
\end{itemize}

\section{Couplings, Maps and Costs}
\label{sec:coup_maps_costs}

In this section we establish some correspondence between couplings and maps $\alpha:\mathbf{T}\rightarrow \mathbf{G}$, giving also the definition of the cost of a coupling.
As a first step, given a coupling $C$ we induce two maps $\alpha_C:V_T-U_C^T\rightarrow \mathbf{G}$ and $\beta_C:V_G-U_C^G\rightarrow \mathbf{T}$, following these rules:

\begin{enumerate}
\item if $(x,y)\in C$, then $\alpha_C(x):=y$; 
\item $x\in D_C^T$ with $\#\Lambda(x)=0$, then
$\alpha_C(x):=s^k_T(\eta(x))$ 
with $k$ being:
	\begin{itemize}
	\item $k=f(x)+\frac{1}{2}\left( f(\varphi(x))-f(x)\right) - g(\eta(x))$ if $g(\eta(x))\leq f(x)+\frac{1}{2}\left( f(\varphi(x))-f(x)\right)$;
	\item $k=0$ otherwise. 
	\end{itemize}
Where the idea is that we want to send $x$ above some coupled vertex changing its height \virgolette{as little as possible}. Bear in mind that the sequence $[x,\varphi(x)]$ should not appear in the image of $\beta\alpha$;
\item if instead $x\in D_C^T$ with $\#\Lambda(x)>1$, then
$\alpha_C(x):=\chi(x)$.
\end{enumerate}

The map $\beta_C$ is obtained in the same way, exchanging the roles of $(T,f)$ and $(G,g)$.
\Cref{fig:coup_to_eps_2} shows an example in which $\beta_C$ is obtained starting from a coupling.

As a first result we obtain that the maps we induce respect the trees structures of the respective merge trees.

\begin{prop}\label{prop:alpha_order}
Consider $x,x'\in V_{T}-U_T$: if $x<x'$ then $\alpha_C(x)\leq \alpha_C(x')$. 
\end{prop}

Now we need to extend $\alpha_C$ to a continuous function between $\mathbf{T}$ and $\mathbf{G}$. To formally extend the map, we need the following lemma.

\begin{lem}\label{lem:sequences}
Any $x\in \mathbf{T}$, $x\notin \pi_T(C)\bigcup D$, is contained in a sequence of edges 
$[l_C(x),u_C(x)]$ such that 
$l_C(x) \in \max \{v\leq x \mid v \in \pi_T(C)\cup D\}$ and $u_C(x)$ is either $u_C(x)=+\infty$ or
$u_C(x) = \min \{v\geq x \mid v \in \pi_T(C)\cup D\}$.

Moreover, $l_C(x)$ is unique if $\{v\in U \mid l_C(x)\leq v \leq x \} = \emptyset$ and, if $\{v\in U \mid l_C(x)\leq v \leq x \} \neq \emptyset$, there exist one and only one $l_C(x)$ such that $l_C(x)\notin D$. 

\end{lem}

Now we can obtain a continuous map $\alpha_C:\mathbf{T}\rightarrow \mathbf{G}$.

\begin{prop}\label{prop:C_extension}
We can extend $\alpha_C$ to a continuous map between $\mathbf{T}$ and $\mathbf{G}$. This map, with an abuse of notation, is still called $\alpha_C$.
\begin{proof}
 
We define $\alpha_C(x)$ for $x\notin \pi_T(C)\bigcup D_C^T$. By \Cref{lem:sequences} we have 
$x\in [l(x),u(x)]$, with $\alpha_C$ being defined for $l(x)$ and $u(x)$. Whenever $l(x)$ is not unique, we take the unique $l(x)\in \pi_T(C)$. Clearly we have $f(l(x))<f(x)<f(u(x))$.
Thus, if $u(x)<+\infty$, there is a unique $\lambda\in [0,1]$ such that $f(x)= \lambda f(l(x)) + (1-\lambda) f(u(x))$. 
Having fixed such $\lambda$, we define $\alpha_C(x)$ to be the unique point in $[\alpha_C(l(x)), \alpha_C(u(x))]$ - which is a sequence of edges thanks to \Cref{prop:alpha_order} - with height $\lambda g(\alpha_C(l(x)))+(1-\lambda)g(\alpha_C(u(x)))$.	
If $u(x)=+\infty$ then $x>v$ with $v=\max \pi_T(C)$. Then we set $\alpha_C(x)=y$ such that $y\geq \alpha_C(v)$ and $g(y)=g(\alpha_C(v))+ f(x)-f(v)$. Note that, for $x,x'$ such that $u(x),u(x')=+\infty$, 
$\alpha_C$ always preserves distances and $g(\alpha_C(x))-f(x)=g(\alpha_C(v))-f(v)$. Thus, on such points it is a continuous function.

Consider now a converging sequence $x_n\rightarrow x$ in $\mathbf{T}$. 
We know that definitively $\{x_n\}_{n\in\mathbb{N}}$ is contained in one or more edges  containing $x$. Thus we can obtain a finite set of converging subsequences by intersecting 
$\{x_n\}$ with such edges. With an abuse of notation from now on we use $\{x_n\}_{n\in\mathbb{N}}$ to indicate any such sequence.
Each of those edges is contained in a unique sequence of edges of the form $[l(x'),u(x')]$, for some $x'$ - up to, eventually, taking $l(x')\notin D$. Thus $\{x_n\}\subset [l(x'),u(x')]$ induces a unique sequence $\{\lambda_n\}\subset [0,1]$ such that $f(x_n)=\lambda_n f(l(x'))+(1-\lambda_n)f(u(x'))$ and $\lambda_n \rightarrow \lambda_x$, with 
$f(x)=\lambda_x f(l(x'))+(1-\lambda_x)f(u(x'))$. 
By \Cref{lem:sequences}, $\alpha_C(x)\in [\alpha_C(l(x')),\alpha_C(u(x'))]$. Moreover, by construction, $g(\alpha_C(x_n))\rightarrow \lambda_x g(\alpha_C(l(x')))+(1-\lambda_x)g(\alpha_C(u(x')))=g(\alpha_C(x))$.
Thus $\alpha_C$ is continuous.

\end{proof}
\end{prop}

In order to start relating maps induced by couplings and $\varepsilon$-good maps, we define the $cost$ of the coupling $C$, which is given in terms of how much $\alpha_C$ moves the points of $\mathbf{T}$.

\begin{defi}
Given $C\in\mathcal{C}(T,G)$ and $x\in \mathbf{T}$, we define $cost_C(x)=\mid g(\alpha_C(x))-f(x)\mid$. 
Coherently we define 
$\parallel C\parallel_\infty = \max\{\parallel g\circ \alpha_C-f\circ \text{Id}_\mathbf{T} \parallel_\infty,\parallel f\circ \beta_C-g\circ \text{Id}_\mathbf{G} \parallel_\infty\}$.
\end{defi}

Note that, given $C\in\mathcal{C}(T,G)$ and $x\in \mathbf{T}$, we have one of the following possibilities:
\begin{itemize}
\item if $(x,y)\in C$, $cost_C(x)=\mid f(x)-g(y)\mid $;
\item if $x\notin\pi_T(C)\bigcup D$, $cost_C(x)\leq \max\{cost_C(l(x)),cost_C(u(x))\}$;
in fact:
\begin{align*}
&\mid g(\alpha_C(x))-f(x)\mid=\\
&\mid \lambda g(\alpha_C(l(x)))-\lambda f(l(x)) + (1-\lambda) g(\alpha_C(u(x)))-(1-\lambda) f(u(x))\leq\\ 
&  \lambda cost_C(l(x)) + (1-\lambda)cost_C(u(x))
\end{align*}
\item we are left with the case $x\in D$. We have two different scenarios:
	\begin{itemize}
	\item if $\#\Lambda(x)=0$, then  $cost_C(x)= \max\{\left( f(\varphi(x))-f(x)\right) /2,g(\eta(x))-f(x)\}$. We point out that $\left( f(\varphi(x))-f(x)\right) /2$ is the cost of deleting the path $[x,\varphi(x)]$ in two steps: we halve the distance between $x$ and $\varphi(x)$ with $\alpha$ and then with $\beta$ we have $f(\beta\alpha(x))=f(\varphi(x))$. This can happen if below $w$ (with $(\delta(x),w)\in C$) there is \virgolette{room} to send $x$ at height $f(x) + \left( f(\varphi(x))-f(x)\right) /2$; if instead $\eta(x)$ is higher than $f(x) + \left( f(\varphi(x))-f(x)\right) /2$ we simply send $x$ to $\eta(x)$; 
	\item if $\#\Lambda(x)>1$, we have $cost_C(x)= \mid f(x)- g(\chi(x))\mid $.
	\end{itemize}
\end{itemize}

As a consequence we point out two facts:
\begin{enumerate}
\item for every $\alpha:\mathbf{T}\rightarrow \mathbf{G}$ such that, for every $x\in \pi_T(C)\bigcup D$, $\alpha(x)=\alpha_C(x)$  we have:
\begin{equation*}
\parallel g\circ \alpha_C-f\circ \text{Id}_\mathbf{T} \parallel_\infty \leq \parallel g\circ \alpha-f\circ \text{Id}_\mathbf{T} \parallel_\infty
\end{equation*}
\item if we fix some total ordering of the elements in $V_T \coprod V_{G}$, so that we can write $V_T \coprod V_{T'}=(a_1,\ldots, a_n)$, then $\parallel C\parallel_\infty = \max(cost(a_1),\ldots, cost(a_n))$.
\end{enumerate}

\section{From Couplings to $\varepsilon$-good Maps}
\label{sec:couplings_maps}

\begin{figure}
    \begin{subfigure}[c]{0.97\textwidth}
    	\centering
    	\includegraphics[width = \textwidth]{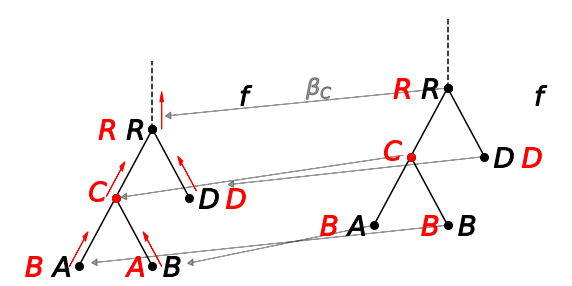}
    	\caption{}
		\label{fig:coup_to_eps_1}	
	\end{subfigure}
	
	    \begin{subfigure}[c]{0.97\textwidth}
    	\centering
    	\includegraphics[width = \textwidth]{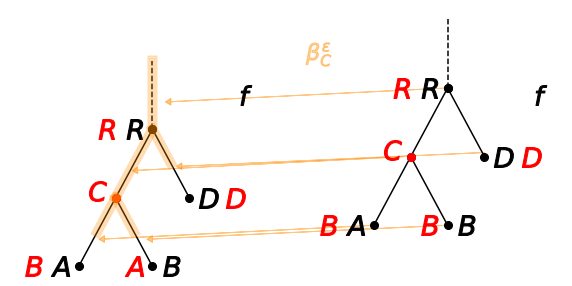}
   	\caption{}
	\label{fig:coup_to_eps_2}	
	\end{subfigure}
\caption{With the same notation/colors used in  \Cref{fig:coup_1} and display the machinery defined in \Cref{sec:couplings_maps}. Start from the subfigure (a). The gray arrows represent the map $\beta_C:V_G\rightarrow V_T$, which then is naturally extended to the metric trees. The image of such map is then shifted upwards using the structure map $s^{k}_T$ as in \Cref{sec:couplings_maps}, to obtain $\beta_C^\varepsilon$, with the shift being indicated by the upwards red arrows. In subfigure (b), the orange arrows - and the shaded orange portion of $\mathbf{T}$ -  represent the $\varepsilon$-good map $\beta_C^\varepsilon$ obtained with the composition of $\beta_C$ with the upward shift.}
\label{fig:coup_to_eps}
\end{figure}

In this section we prove the first fundamental interaction between couplings and the interleaving distance between trees, with our final goal being to replace the problem of finding optimal maps with the combinatorial problem of obtaining optimal couplings.
\Cref{fig:coup_to_eps} summarizes all the steps of the procedure which is formally addressed in this section.
Throughout the section, we fix a coupling $C$ and some $\varepsilon\geq \parallel C\parallel_\infty$. 

We build a map $\alpha^{\varepsilon}_C:\mathbf{T}\rightarrow \mathbf{G}$: $\alpha_C^{\varepsilon}(x)=s^{k_x}_G(\alpha_C(x))$ with $k_x=f(x)+\varepsilon-g(\alpha_C(x))$ depending on $x$.
Now we check that, since $\varepsilon\geq cost(C)$, we always have $k_x\geq 0$. In fact it is enough to check this property on $V_T$:
\begin{enumerate}
\item  for $(x,y)\in C$ we have $\varepsilon \geq \mid f(x)-g(y)\mid $ and thus $f(x)+\varepsilon-g(y)\geq 0$;
\item for $x\in D_C^T$ and  $\#\Lambda(x)=0$, we have $\varepsilon\geq  \max\{(f(\varphi(x))-f(x))/2,g(\eta(x))-f(x) \}$.
Thus $f(x)+\varepsilon \geq g(\eta(x))$ and so $k\geq 0$.
Similarly if $\#\Lambda(x)>1$, we have $w=\chi(x)$ and $k =  f(x) + \varepsilon - g(w) \geq 0$. 
\end{enumerate}

We want to prove that $\alpha^{\varepsilon}_C$ is an $\varepsilon$-good map.
Note that, by construction $g(\alpha_C^{\varepsilon}(x))=g(\alpha_C(x))+f(x)+\varepsilon-g(\alpha_C(x))=f(x)+\varepsilon$.

\begin{rmk}\label{rmk:monotone}
The map $\alpha^{\varepsilon}_C$ is such that, given $x,x'\in \mathbf{T}$ with $x<x'$, we have $\alpha^{\varepsilon}_C(x)<\alpha^{\varepsilon}_C(x')$. 
\end{rmk}

As a first step we obtain that $\alpha^{\varepsilon}_C$ satisfies some necessary conditions in order to be $\varepsilon$-good: it is in fact continuous and moves point \virgolette{upwards} by $\varepsilon$.

\begin{prop}\label{prop:extension}
The map $\alpha^{\varepsilon}_C$ is a continuous map between $\mathbf{T}$ and $\mathbf{G}$ such that for every $x\in\mathbf{T}$ we have $g(\alpha_C^\varepsilon(x))=f(x)+\varepsilon$.

\end{prop}

Before proving the main result of this section we need a short lemma.

\begin{lem}\label{lem:diseq_interl}
Let $(v,w),(v',w')\in C$ and let $x=\LCA(v,v')$, $y=\LCA(w,w')$. Then $\mid f(x)-g(y)\mid \leq  \varepsilon$.

\end{lem}

At this point we are ready to prove the remaining properties which make $\alpha_C^\varepsilon$ an $\varepsilon$-good map.

\begin{teo}\label{teo:eps_to_C}
The map $\alpha_C^\varepsilon$ is an $\varepsilon$-good map.

\end{teo}

\section{From $\varepsilon$-Good Maps to Couplings}
\label{sec:map_couplings}

This time we start from an $\varepsilon$-good map $\alpha:\mathbf{T}\rightarrow \mathbf{G}$ and our aim is to induce a $C$ with $\parallel C \parallel_\infty\leq \varepsilon$.

Given a metric merge tree $\mathbf{T}$ we define the following map: $L:\mathbf{T}\rightarrow V_T$ given by $L(x)= \max\{v \in V_T \mid v\leq x\}$. Note that if $x\in V_T$ then $L(x)=x$, otherwise $x$ is an internal vertex of one edge $(a,b)$, with possibly $b=+\infty$, and $L(x)=a$. 
With this notation we introduce the following maps, which are analogous to $\alpha^{\downarrow}$ and $\beta^{\downarrow}$ defined in \cite{interl_approx_2018} and to other maps defined in the proof of Theorem 1 in \cite{pegoraro2024functional}:

\begin{align}
\phi:&V_T\rightarrow V_{G}\\
&v \mapsto L(\alpha(v))
\end{align}

\begin{align}
\psi:&V_G\rightarrow V_{T}\\
&w \mapsto L(\beta(v))
\end{align}

\begin{figure}
    \begin{subfigure}[c]{0.97\textwidth}
    	\centering
    	\includegraphics[width = \textwidth]{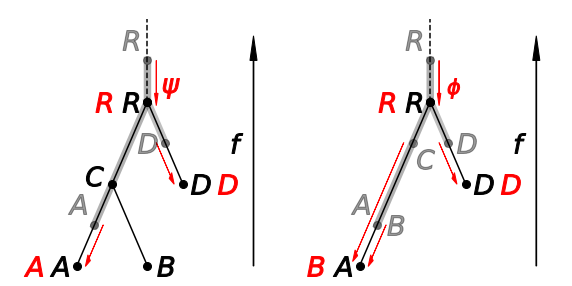}
	\end{subfigure}
\caption{Given the two maps of \Cref{fig:interleaving}, the shaded gray represent the image of those maps, the red arrows give the maps $\psi$ (left) and $\phi$ (right) - see \Cref{sec:map_couplings}, and the red letters next to the black ones indicate two possible couplings (one in leftmost tree and one in the rightmost tree) satisfying \Cref{teo:C_to_eps}. These couplings are compatible with the procedure outlined  in the proof of \Cref{teo:C_to_eps} upon perturbing the leftmost tree to meet the generality condition (G).}
\label{fig:alpha_beta_coupling}
\end{figure}

Note that we always have $g(\phi(v))\leq g(\alpha(v))\leq f(v)+\varepsilon$.
These maps will be keys in the proof of the upcoming theorem since they will help us in closing the gap between the continuous formulation of the interleaving distance and the discrete matching of merge trees via couplings.
The reader should refer to \Cref{fig:alpha_beta_coupling}
for a visual example of the above definitions.

We prove a corollary which characterizes the maps we just defined, and will be used in what follows.

\begin{cor}\label{lem:lemma_0}
Let $v,v'\in V_T$; if $v<v'$ then $\phi(v)\leq \phi(v')$. 
\end{cor}

Clearly an analogous result holds for $\psi$, which implies that, in the setting of the corollary, $\psi(\phi(v))\leq \psi(\phi(v'))$.

Now we prove the main result of this section.

\begin{teo}	\label{teo:C_to_eps}
Given $\alpha$ (and $\beta$) $\varepsilon$-good maps between $\mathbf{T}$ and $\mathbf{G}$, then there is a coupling $C$ between $T$ and $G$ such that $\parallel C \parallel_\infty\leq \varepsilon$.
\end{teo}

Putting together \Cref{teo:eps_to_C} and \Cref{teo:C_to_eps}, we see that two merge trees are $\varepsilon$-interleaved if and only if there is a coupling $C$ between the tree such that $\parallel C \parallel_\infty = \varepsilon$. 
Thus, computing the interleaving distance amounts to finding a minimal-cost coupling between two merge trees. As a byproduct of this result, we obtain another proof that the interleaving distance between metric merge tree can be found as a minimum and not just as an infimum as in \Cref{def:interleaving}, for the set of available couplings is finite.

We close this section with a claim which we want to investigate in future works, which stems from one of the definitions of the Gromov-Hausdorff distance \cite{adams2022gromov} and the already exploited link between such distance in the case of metric merge trees and the interleaving distance \citep{interl_approx}.
Given $\alpha:\mathbf{T}\rightarrow \mathbf{G}$ and $\mu:\mathbf{T}\rightarrow \mathbf{T}$ if we define
$H(\alpha):=\max_{x\in\mathbf{T}}\mid g(\alpha(x))-f(x)\mid$ and $2\cdot D(\mu)= \max_{x\in\mathbf{T}}d(x,\mu(x))$ (where $d$ is the shortest path metric on $\mathbf{T}$) we claim that:
\[
d_I(\mathbf{T},\mathbf{G}) = \min_{\alpha,\beta} \{H(\alpha),H(\beta), D(\beta\alpha),D(\alpha\beta)\}
\] 

where $\alpha:\mathbf{T}\rightarrow \mathbf{G}$ and $\beta:\mathbf{G}\rightarrow \mathbf{T}$ are continuous monotone maps. 

\section{Properties of Couplings}
\label{sec:prop}
In this section we present two different properties of the interleaving distance which are useful for the approximation scheme we propose in \Cref{sec:LIP}.

\subsection{Decomposition Properties}
\label{sec:decomposition_interl}

We exploit the equivalent formulation of the interleaving distance via couplings to prove that it can be approximated
via the solution of in smaller independent subproblems which are then aggregated to solve the global one. 

The main result of the section involves couplings with some strong properties, which we call \emph{special} and collect under the following notation:

\begin{equation*}
\mathcal{C}^o(T,G):=\{C\in \mathcal{C}(T,G)\mid \arg\min_{v<\max \pi_T(C)} f(v) \in \pi_T(C)\text{ and }\arg\min_{w<\max \pi_G(C)} g(w) \in \pi_G(C)\}.
\end{equation*}

Before proceeding, we need few pieces of notation used to lighten the dissertation and one last technical definition.
Given $x\in V_T$ and $y\in V_G$ we define $T_x=sub_T(x)$ and $G_y=sub_G(y)$. 
Moreover 
$\mathcal{C}_R(T,G):=\lbrace C\in \mathcal{C}(T,G) \mid (r_T,r_G) \in C\rbrace$. Similarly, we have $\mathcal{C}^o_R(T,G):=\mathcal{C}_R(T,G)\bigcap \mathcal{C}^o(T,G)$.

\begin{defi}
Given a partially ordered set $(A,<)$, $A$ is antichain if for any $a,b\in A$, $a\neq b$, there is not an element $c$ such that $c>a$ and $c>b$.
\end{defi}

\begin{figure}
    \begin{subfigure}[c]{0.97\textwidth}
    	\centering
    	\includegraphics[width = \textwidth]{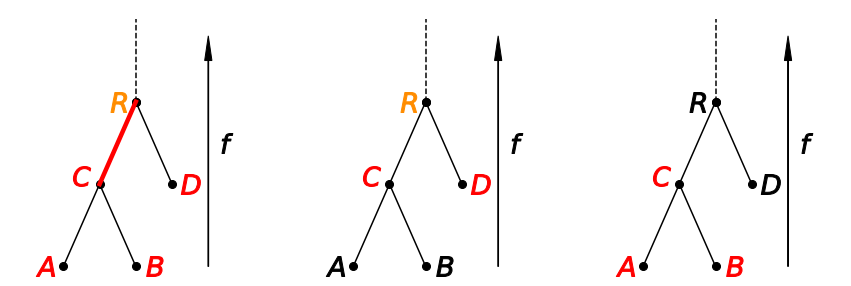}
	\end{subfigure}
\caption{The letters in red show three examples of $C^*\subset V_T\times V_G$ satisfying property $(A0)$ - see \Cref{defi:antichain}. However the rightmost example does not satisfy the antichain condition $(A1)$  - see \Cref{defi:antichain}, for there is $C>A$ and $C>B$. In orange we highlight the vertices giving $\LCA(\pi_T(C^*))$. The green branch in the leftmost tree signifies that, since $\#\Lambda_{C^*}(C)=2$, $C$ is deleted. }
\label{fig:antichain}
\end{figure}

Now we introduce the combinatorial objects we will use to decompose the following optimization problem: $\min_{C\in\mathcal{C}(T,G)} \parallel C\parallel_\infty$. The reader may look at \Cref{fig:antichain} to find examples involving the following definition.

\begin{defi}\label{defi:antichain}
We define $\mathcal{C}^*(T,G)$ as the set of $C^*\subset V_T\times V_G$ such that:
\begin{itemize}
\item[(A0)] $C^*\bigcup 
\big \lbrace \left( \LCA(\pi_T(C^*)),\LCA(\pi_G(C^*))
 \right) \big \rbrace \in \mathcal{C}(T,G)$;
\item[(A1)] $C^*$ is an antichain.
\end{itemize}  
\end{defi}

The idea behind the definition of $\mathcal{C}^*(T,G)$
is that we are assuming that we already know $\arg\min \{\parallel C'\parallel_\infty \mid C'\in\mathcal{C}_R(T_x,G_y)\}$ for every $x\in E_T$ and $y\in E_G$, and we use $C^*\in \mathcal{C}^*(T,G)$ to optimally aggregate these results to find the optimal global coupling between $T$ and $G$.
We formalize such concepts in the following definition.

\begin{defi}
Let $C_{x,y}\in \mathcal{C}_R(T_x,G_y) $ for every $x,y\in V_T\times V_G$. Given $C^*\in\mathcal{C}^*(T,G)$, with $r=\LCA(\pi_T(C^*))$ and $r'=\LCA(\pi_G(C^*))$, we define the \emph{extension} of $C^*$ by means of $\{C_{x,y}\}_{(x,y)\in C^*}$ as $e(C^*)= \{(r,r')\}\bigcup_{(x,y)\in C^*} C_{x,y}$. The set of all possible extensions of $C^*$ is called $E(C^*)$.
If $C_{x,y}\in \arg\min \{\parallel C'\parallel_\infty \mid C'\in\mathcal{C}_R(T_x,G_y)\}$ for all $x,y$ then we call the extension minimal. We collect all minimal extensions of $C^*$ in the set $E_m(C^*)$ - which is never empty. 
If all  $C_{x,y}$ and $e(C^*)$ are special couplings then we call the extension special. We collect all special extensions of $C^*$ in the set $E^o(C^*)$ - which is never empty.
We name $E^o_m(C^*)=E_m(C^*)\cap E^o(C^*)$ the set of minimal special extensions of $C^*$ - this set could be empty.
\end{defi}

Since extensions are couplings, it is obvious that:
\begin{equation}\label{eq:approx_min}
d_I(T,G)\leq \min_{C^*\in\mathcal{C}^*(T,G)}\min_{C \in E_m(C^*)} \parallel C \parallel_\infty.
\end{equation}
\begin{equation}\label{eq:approx_sp}
d_I(T,G)\leq \min_{C^*\in\mathcal{C}^*(T,G)}\min_{C \in E^o(C^*)} \parallel C \parallel_\infty.
\end{equation}

\begin{figure}
    \begin{subfigure}[c]{0.97\textwidth}
    	\centering
    	\includegraphics[width = \textwidth]{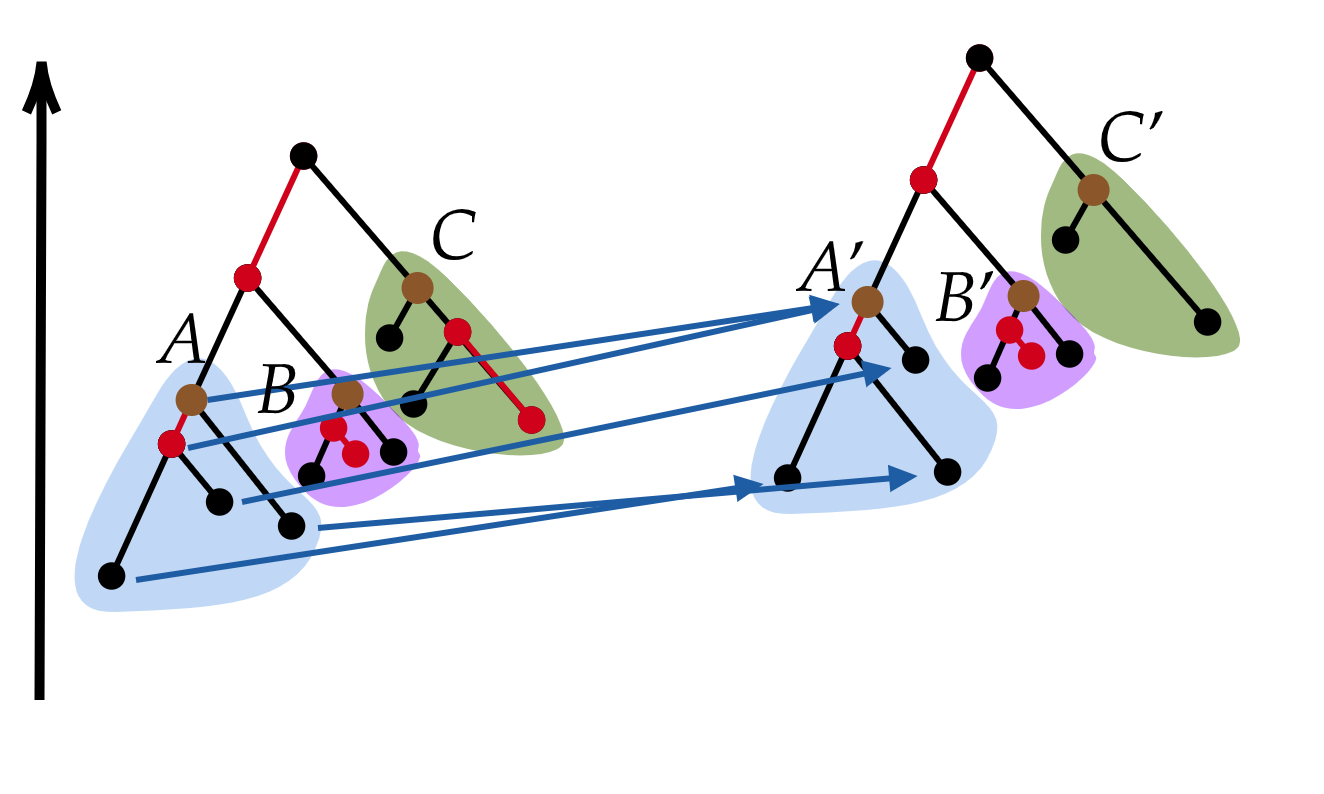}
    	\caption{The couples $\{(A,A'),(B,B'),(C,C')\}$ form antichains via $\pi_T$ and $\pi_G$. The shaded regions indicate how the subtrees rooted in the vertices of the antichain are matched. Minimal couplings are displayed with deletions in red. The remaining vertices can be coupled by visual inspection. We have highlighted the couples of the blue subtree just as an example. Putting together all the couples we obtain an extension of the antichain.}
		\label{fig:decomposition_interl}	
	\end{subfigure}
	
	    \begin{subfigure}[c]{0.97\textwidth}
    	\centering
    	\includegraphics[width = \textwidth]{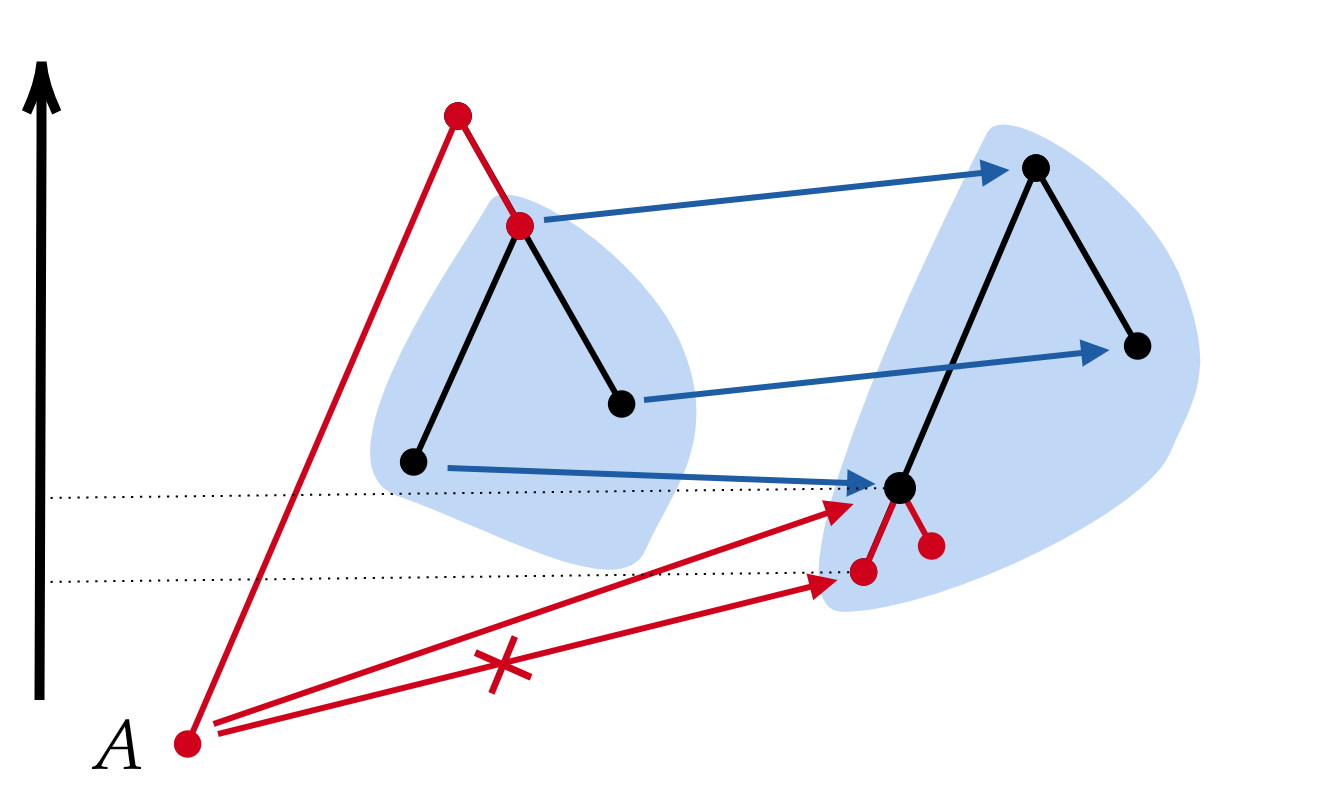}
   	\caption{In this figure we can see a non special extension as the lower vertex of the tree on the right is not matched. Even though the coupling displayed by the blue arrows is special for the blue subtrees, it does not extend to a special coupling between the two trees: the deletion of $A$ should go through the lowest vertex of the right tree, while it is forced to go though its father. As signified by the red arrows. We can also appreciate that minimal extensions are not optimal couplings in general.}
	\label{fig:bad_extension}	
	\end{subfigure}
\caption{Two figures related to decomposition and extensions of couplings.}
\label{fig:extension}
\end{figure}

See also \Cref{fig:decomposition_interl}. Moreover it is also clear that any coupling is an extension of some $C^*\in\mathcal{C}^*(T,G)$ and thus:
\begin{equation*}
d_I(T,G) = \min_{C^*\in\mathcal{C}^*(T,G)}\min_{C \in E(C^*)} \parallel C \parallel_\infty.
\end{equation*}

The upcoming theorem states that there strong relationships between $d_I$ and extensions obtained via a fixed family of $C_{x,y}\in\mathcal{C}_R(T_x,G_y)$.

\begin{teo}[Decomposition]\label{teo:deco}
Consider two merge trees $T$ and $G$ and take
 a collection $\{C_{x,y}\}_{(x,y)\in V_T \times V_G}$ with $C_{x,y}\in \arg\min \{\parallel C'\parallel_\infty \mid C'\in\mathcal{C}_R(T_x,G_y)\}$.
Given $C^*\in\mathcal{C}^*(T,G)$ we name $e(C^*)$ the extension obtained by means of $\{C_{x,y}\}$.
We have:
\begin{equation}\label{eq:low_bound}
\min_{C^*\in\mathcal{C}^*(T,G)} \max \{ cost_{e(C^*)}(v)\mid v \in \pi_T(e(C^*))\text{ or } \#\Lambda_{e(C^*)}(v)>0\}\leq d_I(T,G)
\end{equation}
Moreover, if $C_{x,y}$ is also a special coupling for every $x,y$, we have:
\begin{equation}\label{eq:equal_interl}
d_I(T,G)= \min_{C^*\in\mathcal{C}^*(T,G)}
\parallel e(C^*) \parallel_\infty.
\end{equation}
\end{teo}

We remark that \Cref{teo:deco} is in some sense unexpected: if \Cref{eq:approx_min} and \Cref{eq:approx_sp} are in some sense trivial, \Cref{eq:low_bound} and \Cref{eq:equal_interl} certainly are not. We highlight that the couplings $\{C_{x,y}\}$ are independently fixed at the beginning and not chosen with some joint optimization strategy. Moreover, the proof  of \Cref{eq:equal_interl} depends strongly on the possibility to find optimal couplings which are also special.

\subsection{Approximation by Pruning Operators}
\label{sec:pruning_interl}

Now we prove a second property of the interleaving distance which is very useful when looking for approximations of $d_I$. In particular, we find a way to approximate the distance sensibly reducing the computational complexity of the problem by removing leaves and computing distances between smaller trees.

We briefly introduce the pruning operator $P_\varepsilon$.
Given a merge tree $T$ and $\varepsilon>0$, the merge tree $P_\varepsilon(T)$ is obtained through  a recursive procedure: given a leaf $x$ and its father $x'$, if $f(x')-f(x)<\varepsilon$ we say that $x$ is a small-weight leaf; we want to remove all small-weight leaves - and their fathers if they become order $2$ vertices - from $T$ unless two or more of them are siblings, i.e. children of the same father. In this case we want to remove all leaves but the one being the lowest leaf.
To make this procedure well defined and to make sure that, in the end, no small-weight leaves are left in the tree, we need to choose some ordering of the leaves and to resort to recursion. 
\begin{itemize}
\item[(P)] Take a leaf $l$ such that $f(father(l))-f(l)$ is minimal; if $f(father(l))-f(l)<\varepsilon$, remove $l$ and its father if it becomes an order $2$ vertex after removing $l$. 
\end{itemize}

\begin{figure}
	\centering
    \begin{subfigure}[c]{0.47\textwidth}
    	\centering
    	\includegraphics[width = \textwidth]{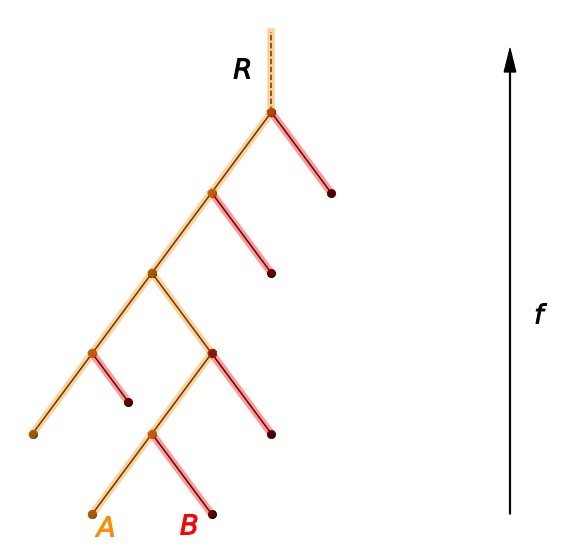}
	\end{subfigure}
\caption{An example of a pruning operator applied on a merge tree $T$. The red branches are removed from the tree, while the orange ones are kept and represent the metric merge tree $P_\varepsilon(T)$. We highlight that the root of the tree changes. The superimposition of the orange shaded tree on the black one gives the injection between them. Lastly, note that deleting $A$ instead of $B$ (these are small-weight siblings with same weights - violating (G)) would give isomorphic merge trees.}
\label{fig:pruning_interl}
\end{figure}

We take $T_0=T$ and apply operation (P) to obtain $T_1$. On the result we apply again (P) obtaining $T_2$ and we go on until we reach a fixed point $T_n = T_{n+1}=P_\varepsilon(T)$ of such sequence. To shed some light on this definition we prove the following lemma; \Cref{fig:pruning_interl} can be helpful in following the statements.

\begin{lem}
\label{lem:pruning}
Given $T$, $\varepsilon>0$ and $P_\varepsilon(T)$, we have a natural injective map $V_{P_\varepsilon(T)} \hookrightarrow V_T$ which identifies vertices in  $P_\varepsilon(T)$ with vertices in $T$. Call $C_\varepsilon\subset V_{P_\varepsilon(T)} \times V_T$ the set of couples induced by $V_{P_\varepsilon(T)} \hookrightarrow V_T$. 
 The following hold:
\begin{enumerate}
\item $C_\varepsilon$ is a coupling;
\item $L_{P_\varepsilon(T)}\subset L_T$ and
for every $v$ and $v'$ such that $v'<v$ and $f(v)-f(v')\geq \varepsilon $, there is $l \in L_{P_\varepsilon(T)}$ such that $\LCA(l,v')< v$;
 in particular $\arg\min f \in L_{P_\varepsilon(T)}$;
\item for every $v\in V_T- V_{P_\varepsilon(T)}$ we have $\#\Lambda_{C_\varepsilon}(v)\leq 1$; in particular if $v\in D$, we have $\#\Lambda(v)=0$ and $f(\varphi_{C_\varepsilon}(v))-f(v)< \varepsilon$;
\item the map:
$\eta_{C_\varepsilon}:D^T_{C_\varepsilon}\rightarrow V_T$ such that $f(\eta_{C_\varepsilon}(x))<f(x)$ for all $x\in D^T_{C_\varepsilon}$;
\item $\parallel C_\varepsilon \parallel_\infty \leq \varepsilon/2$.
\end{enumerate}
\end{lem}

Having characterized the pruned trees with the above lemma, we can obtain the following proposition.

\begin{prop}
\label{prop:pruning}
Given two merge trees $T$ and $G$, we have:
\[
 d_I(T,G)\leq \max\{d_I(P_\varepsilon(T),P_\varepsilon(G)),\varepsilon/2\}
\]

\end{prop}

As a result, if the number of leaves of $T$ and $G$ is too high, we know that we can prune them, reducing the computational complexity of $d_I$, to obtain an estimate from above 
of $d_I(T,G)$. 

We close this section with a claim we would like to investigate in the future.

\begin{claim}
The coupling $C_\varepsilon$ is a minimizing coupling.
\end{claim}


\section{Approximating the Interleaving Distance: Linear Binary Optimization}
\label{sec:LIP}

In this section we exploit \Cref{teo:deco} to obtain 
a dynamical approach to approximate an optimal coupling between two merge trees, by solving recursively linear binary programming problems.

\subsection{Computing the Cost of a Coupling Extension}

As a first step we separate the problem of finding $C^*\in\mathcal{C}^*(T,G)$ with a cost-minimizing extension into two separated problems: the first one is to find a minimal-fixed root coupling, the second one is to compute the cost of deleting the vertices which are not in the subtrees given by the roots fixed by $C^*$. To make the upcoming lemma more clear we establish the following notation: if $r=\LCA(\pi_T(C^*))$ and $r'=\LCA(\pi_G(C^*))$, then for $v\nleq r$, $v_r= \min\{v'\geq r \text{ and } v' \geq v\}$ and $f_r = \min_{v\leq r}f(v)$. Lastly: 
\[H_{r,r'}=\max\{\max_{v\nleq r}0.5\cdot(f(v_r)-f(v)),
\max_{v\nleq r}g_{r'}-f(v),
\max_{w\nleq r'}0.5\cdot(g(w_{r'})-g(w)),
\max_{w\nleq r'}f_r-g(w)\}.
\]

Note that $H_{r,r'}$ does not depend on the chosen extension of $C^*$, and, in fact, it depends only on $r,r'$.  The upcoming lemma states that $H_{r,r'}$  accounts for the deletions of all the vertices which are not below $r$ or $r'$. 

\begin{lem}
\label{lem:extension_cost}
Consider $T,G$ merge trees and $C^*\in\mathcal{C}^*(T,G)$, with $r=\LCA(\pi_T(C^*))$ and $r'=\LCA(\pi_G(C^*))$. Let $C_{r,r'}$ be any extension of $C^*$ and let $C_{r,r'}'$ be a coupling in $\mathcal{C}(T_r,G_{r'})$ made by the same couples as $C_{r,r'}$. Then:
\begin{equation}
\parallel C_{r,r'} \parallel_\infty \geq \max\{\parallel C'_{r,r'} \parallel_\infty, H_{r,r'}\}. 
\end{equation}  If $C_{r,r´}$ is a special coupling, the inequality becomes and equality.
\end{lem}

Since the computation of $H_{r,r'}$ can be easily accessed in a greedy fashion, from now on, we focus on find an approximate solution for $C\in \arg\min \{\parallel C'\parallel_\infty \mid C'\in\mathcal{C}_R(T,G)\} $.

\subsection{Iterative Approach} 

Now we start working out a procedure to approximate $C\in \arg\min \{\parallel C'\parallel_\infty \mid C'\in\mathcal{C}_R(T,G)\} $.
The assumption of our approach is that we already have computed the couplings 
$C_{x,y}\in \mathcal{C}_R(T_x,G_y)$ that we want to use to extend any $C^*\in\mathcal{C}^*(T,G)$, as in \Cref{teo:deco}.
So for instance, if we want to work with special extensions, we assume that we have
$C_{x,y}\in \arg\min\{ \parallel C'\parallel_\infty \mid C'\in \mathcal{C}_R^o(T_x,G_y)\}$ for all $x\in E_T$ and $y\in E_G$ and we exploit them to obtain a special extension of some $C^*\in\mathcal{C}^*(T,G)$.
If we want to work with minimal extensions, instead, we assume to have $C_{x,y}\in \arg\min\{ \parallel C'\parallel_\infty \mid C'\in \mathcal{C}_R(T_x,G_y)\}$. 
We anticipate that both kinds of extensions are important for our purposes as special extension will be used in the following to produce upper bounds for the interleaving distance, while minimal extension will be used accordingly to \Cref{teo:deco} to obtain a lower bound.

We fix the following notation: $f_x = \min_{v\leq x} f(v) $
and $g_y = \min_{w\leq y} g(w) $.
Lastly, we fix a constant $K>0$ such that $K> \max_{x\in V_T,y\in V_G} \mid f(x)-g(y) \mid$. 
In \Cref{sec:approx} we point out which steps of the upcoming procedure may produce errors.
 
\subsection{Variables and Constraints}
\label{sec:constr_interl}
We consider the following set of binary variables: $a_{x,y}$ for $x\in E_T$ and $y\in E_G$; $u_x$ for $x\in E_T$ and $u_y$ for $y\in E_G$. 
The vector of all variables (upon choosing some ordering) will be referred to as $\mathcal{V}$.

We also introduce the following auxiliary variables, which are linear functions of $a_{x,y}$, $u_x$ and $u_y$:
\begin{itemize}
\item $c_x= \sum_y a_{x,y}$ and $c_y= \sum_x a_{x,y}$
\item $\Lambda_x = \sum_{v\leq x} c_v$ and $\Lambda_y = \sum_{w\leq y} c_w$
\item $d_x = \sum_{v\leq x} c_v + \sum_{v> x} c_v$ and $d_y = \sum_{w\leq y} c_w + \sum_{w> y} c_w$
\end{itemize}
and the following linear constraints:
\begin{itemize}
\item[(1)] for every $l\in L_T$: $\sum_{l\leq x < r_T} c_x\leq 1$ and for every $l'\in L_G$: $\sum_{l'\leq y < r_G} c_y\leq 1$
\item[(2)] $u_x \leq 0.5\Lambda_x$ and $u_y \leq 0.5\Lambda_y$ for every $x$ and $y$;
\item[(3)] $u_x \geq m_x\Lambda_x +q_x$ and $u_y \geq m_y\Lambda_y+q_y$ for every $x$ and $y$. With $m_x,q_x$ being any pair of $(m,q)$ such that the following are satisfied: $q<0$, $m+q<0$, $2m+q>0$, $m(\#L_T)<1$ (analogously for $m_y,q_y$);
\item[(4)] (only for special extensions) let $\widetilde{x}=\arg \min_{v\in V_T} f(v)$ and $\widetilde{y}=\arg \min_{w\in V_G} g(w)$; then we ask $\sum_{\widetilde{x}\leq x < r_T} c_x\geq 1$ and $\sum_{\widetilde{y}\leq y < r_G} c_y\geq 1$.
\end{itemize}

The set of vectors of variables which satisfy these constraints is called $\mathcal{K}_m(T,G)$ or $\mathcal{K}^o(T,G)$ depending on whether (4) are, respectively, included or not.
Note that $\mathcal{K}_m(T,G)\supset \mathcal{K}^o(T,G)$. To lighten the notation, we often avoid explicit reference to $T$ and $G$ when talking about the set of possible solutions, and, when we do not wish to distinguish between $\mathcal{K}_m$ and $\mathcal{K}^o$, we simply write $\mathcal{K}$.
Now we briefly comment on the variables and constraints to try to give some intuition about their roles:

\begin{itemize}
\item the variables $a_{x,y}$ are used to match $x$ with $y$, i.e. add the couple $(x,y)$ to the coupling. In particular constraint (1) ensures that if $\mathcal{V}\in\mathcal{K}$, the set $C(\mathcal{V}):=\{(x,y)\mid v_{x,y}=1\}$ 
belongs to $\mathcal{C}^*(T,G)$;
\item the variables $d_x$ and $u_x$ are used to determine if $x$ must be deleted. Start with $d_x$. The main idea which our optimization procedure is based on, is that if $a_{v,w}=1$ with $v>x$, it means that $x$ is taken care of by the coupling $C_{v,w}$ and thus we want to \virgolette{ignore} such $x$. In other words, having $d_x=0$, implies that $x$ is deleted with $\#\Lambda(x)=0$;
\item now we turn to $u_x$. Observe that, if $\Lambda_x = 0,1$, then $u_x=0$ while if $\Lambda_x\geq 2$, $u_x = 1$. Note that $\Lambda_x\leq \#L_T$ and $\Lambda_y\leq \#L_G$. Constraints (2) and (3) are fundamental to fix the value of $u_x$, depending linearly on $\Lambda_x$: $u_x=1$ if and only if $x\geq \LCA(v,v')$, with $a_{v,w}=1$ and $a_{v',w'}=1$, for some $w,w'\in V_G$. 
Thus having $u_x=1$ means that $x$ is deleted with $\#\Lambda(x)>1$;
\item lastly, we comment on constraints (4). These constraints are  asking that there is one point above the lowest vertex in each tree which is coupled by $C(\mathcal{V})$; this in particular implies that, if we assume  $C_{x,y}\in\mathcal{C}^o(T_x,G_y)$, then $\widetilde{x}$ and $\widetilde{y}$ (with $\widetilde{x}=\arg \min_{v\in V_T} f(v)$ and $\widetilde{y}=\arg \min_{w\in V_G} g(w)$) are never deleted. Thus,  $\{(r_T,r_G)\}\bigcup_{(x,y)\in C(\mathcal{V})} C_{x,y}$ is a special extension, i.e. a coupling in $\mathcal{C}_R^o(T,G)$.  
\end{itemize}

As a consequence of these observations, any $C\in \mathcal{C}^*(T,G)$ 
induces a unique vector of variables $\mathcal{V}\in\mathcal{K}$, such that $C(\mathcal{V})=C$. Moreover if $C_{x,y}\in\mathcal{C}^o(T_x,G_y)$ and $\mathcal{V}\in\mathcal{K}^o$, then $\{(r_T,r_G)\}\bigcup_{(x,y)\in C(\mathcal{V})} C_{x,y}$ is a special extension of $C(\mathcal{V})$.

\subsection{Objective Functions}
\label{sec:obj_fn}

A key fact that we highlight is that, by property (A1), if $x=\LCA(v,v')$, with $a_{v,w}=1$, $a_{v',w'}=1$ and $x<r_T$, then $\delta(x)=r_T$ due to the antichain condition. So we know that any vertex $x'$ such that $x\leq x' <r_T$, is in $D_{C(\mathcal{V})}^T$. 
Thus, given $x\in V_T$, and with $x_f$ being its father, the \virgolette{local} objective function, depends on the following quantities:

\begin{itemize}
\item $\sum_{y} a_{x,y}\parallel C_{x,y}\parallel _\infty$: it is the cost of matching $T_x$ and $G_y$, if $(x,y)$ is added to $C$. If $C_{x,y}$ is not a minimal coupling, this is an upper bound; 
\item $\mid g(r_G)-f(x)\mid u_x$: it is an upper bound to the cost of deleting $x$ with $\#\Lambda(x)>1$; 
\item $A_x=0.5(f(x_f)-f_x)(1-d_x)$: in case $x$ is deleted with $\#\Lambda(x)=0$, it gives a lower bound to the deletion of $T_x$, as it is equal to deleting the lowest point below $x$ to the height of the father of $x$, in two step. It is an exact value if $x_f=\varphi(x)$  and $g(\eta(x))-f_x< A_x$;
\item $B_x = \lbrace \sum_y a_{v,y}\left( g_y-f_x \right) - K d_x\mid v<x_f \rbrace$: in case $x$ is deleted with $\#\Lambda(x)=0$ it helps giving an upper bound to the deletion of $T_x$, depending on what happens below $x_f$. If there is $v<x_f$ such that $v$ is coupled with $y$ and $C_{x,y}$ is a special coupling, then we know that $g_y\geq g(\eta(x))$, and so $g_y-f_x \geq g(\eta(x))-f_x$. Thus,  if $C_{x,y}$ is a special coupling, $\max \{\max B_x,A_x\}\geq cost(x)$. 
\end{itemize}

For any $x\in E_T$ we define: 
\[
\Gamma^{\uparrow}(x)=\max \left( \sum_{y} a_{x,y}\parallel C_{x,y}\parallel _\infty, (g(r_G)-f(x))u_x,
A_x,\max B_x \right)
\]
and:
\[
\Gamma^{\uparrow}(y)=\max \left((f(r_T)-g(y))u_y,
A_y,\max B_y \right)
\]
for any $y\in E_G$.

In full analogy we set:
\[
\Gamma^{\downarrow}(x)=\max \left( \sum_{y} a_{x,y}\parallel C_{x,y}\parallel _\infty, 
A_x\right)
\]
and:
\[
\Gamma^{\downarrow}(y)=
A_y
\]
for any $y\in E_G$.

Lastly: $\Gamma^{\uparrow}(r_T)=\Gamma^{\downarrow}(r_T)=\mid r_T-r_G\mid $.

With an abuse of notation we write:

\begin{align}\label{eq:cost}
&\Gamma^{\uparrow}(\mathcal{V}) = \max_{x\in V_T\coprod V_G}\Gamma^{\uparrow}(x)\\
&\Gamma^{\downarrow}(\mathcal{V}) = \max_{x\in V_T\coprod V_G}\Gamma^{\downarrow}(x).
\end{align}

Note that all the functions we defined in this section depend on $T$, $G$ and  on the family $\{C_{x,y}\}_{(x,y)\in V_T\times V_G}$ but we avoid explicit reference to those dependencies for notational convenience.

We sum up the main properties of the definitions we have just stated with the following proposition.

\begin{prop}\label{prop:opt}
Given $C\in \mathcal{C}^*(T,G)$ and $C_{x,y}\in\mathcal{C}_R^o(T_x,G_y)$ for all $(x,y)\in C$, we call 
\begin{equation*}
C^o:=\{(r_T,r_G)\} \bigcup_{(x,y)\in C} C_{x,y}.
\end{equation*}
If $C^o$ is a special extension then:
\begin{equation}
\parallel C^o\parallel_\infty\leq \Gamma^{\uparrow}(\mathcal{V})
\end{equation} 

with $\mathcal{V}$ being the unique vector in $\mathcal{K}^o$ such that $C(\mathcal{V})=C^o$.

Viceversa, given $C\in \mathcal{C}^*(T,G)$, and $C_{x,y}\in\mathcal{C}_R(T_x,G_y)$ with minimal cost for all $(x,y)\in C$, we call 
\begin{equation*}
C_m:=\{(r_T,r_G)\} \bigcup_{(x,y)\in C} C_{x,y}.
\end{equation*}
Then:
\begin{equation}
\Gamma^{\downarrow}(\mathcal{V}) \leq \max \{ cost_{C_m}(v)\mid v \in \pi_T(C_m)\text{ or } \#\Lambda_{C_m}(v)>0\}
\end{equation} 

with $\mathcal{V}$ being the unique vector  in $\mathcal{K}_m$ such that $C(\mathcal{V})=C_m$.
\end{prop}

Putting together \Cref{teo:deco} and \Cref{prop:opt} we obtain as a corollary:

\begin{cor}
\label{cor:approx_1}
Consider $T,G$ merge trees, and take:
\begin{enumerate}
    \item a collection of $C_{x,y}\in \arg\min \{\parallel C'\parallel_\infty \mid C'\in\mathcal{C}_R(T_x,G_y)\}$.
    \item a collection of $C'_{x,y}\in \arg\min \{\parallel C'\parallel_\infty \mid C'\in\mathcal{C}^o_R(T_x,G_y)\}$.
\end{enumerate}
Then:

\begin{equation*}
\min_{\mathcal{V}\in \mathcal{K}_m}\Gamma^{\downarrow}(\mathcal{V})\leq \min \{\parallel C \parallel_\infty \mid C\in\mathcal{C}_R(T,G)\}\leq \min_{\mathcal{V}\in \mathcal{K}^o}\Gamma^{\uparrow}(\mathcal{V})
\end{equation*}
where $\Gamma^{\downarrow}$ is computed with the collection $\{C_{x,y}\}$ and $\Gamma^{\uparrow}$ is computed with the collection $\{C'_{x,y}\}$.
\end{cor}

Now we get rid of the fixed roots, obtaining an approximation for $d_I(T,G)$ by putting together \Cref{teo:deco}, \Cref{lem:extension_cost} and \Cref{prop:opt}.

\begin{cor}\label{cor:opt}
In the same setting as \Cref{cor:approx_1}, for each $(x,y)\in V_T\times V_G$ we have:
\begin{equation*}
W_{x,y}^\downarrow:=\min_{\mathcal{V}\in \mathcal{K}_m}\Gamma^{\downarrow}(\mathcal{V})\leq \min \{\parallel C\parallel_\infty \mid C\in\mathcal{C}_R(T_x,G_y)\}\leq W_{x,y}^\uparrow:=\min_{\mathcal{V}\in \mathcal{K}^o}\Gamma^{\uparrow}(\mathcal{V}),
\end{equation*}
with the constraints defining $\mathcal{K}_m$ and $\mathcal{K}^o$ depending on the subtrees $T_x$ and $G_y$. Moreover $\Gamma^{\downarrow}$ is computed with $\{C_{x,y}\}$ and $\Gamma^{\uparrow}$ is computed with $\{C'_{x,y}\}$.

Consequently:
\begin{equation}\label{eq:bounds}
\min_{(x,y)\in V_T\times V_G}
\max\{H_{x,y},W^\downarrow_{x,y}\}
\leq d_I(T,G)
\leq 
\min_{(x,y)\in V_T\times V_G}
\max\{H_{x,y},W^\uparrow_{x,y}\}.
\end{equation}
\end{cor}

\subsection{Approximation Errors} 
\label{sec:approx}

We take this section to briefly isolate which are the situations in which our procedure may produce errors w.r.t. the true interleaving distance.

\begin{enumerate}
\item $\mid g(r_G)-f(x)\mid u_x$: clearly $r_G$ (resp. $r_T$) is an upper bound for $\chi(x)$ (resp. $\chi(y)$) with $x$ (resp. $y$) being deleted with $\#\Lambda(x)>1$ (resp. $\#\Lambda(x)>1$) and so  $\mid g(r_G)-f(x)\mid u_x$ is an upper bound to the cost of the corresponding deletion.
\item $B_{x'}$: for $v\in D_T$ with $\#\Lambda(v)=0$, $x=\varphi(v)$ and $x=father(x')$ we may have $\mid f(v)- g(\eta(v))\mid \leq \max B_{x'}$.
\end{enumerate}

We also point out the following fact.
In the definition of $\Gamma^\uparrow$ the biggest potential source of error are the terms $\mid g(r_G)-f(x)\mid u_x$ which make very costly the deletion of internal vertices in order to swap father children relationships. This is an issue which, to some extent, brings together our approximation scheme for the upper bound of the interleaving distance and other distances for merge trees which have been defined in literature, such as 
\cite{merge_farlocca, merge_farlocca_2, merge_wass}. A very detailed explanation of the problems arising from this fact can be found in \cite{pegoraro2024finitely}, along with the solutions that the authors propose to mitigate them. All these solutions apply also to our case, but, on top of them, we have a couple of additional advantages:
1) in the definition of $\Gamma^\uparrow$ we could also replace $\mid g(r_G)-f(x)\mid u_x$ with $\mid f(father(x))-f(x)\mid u_x$. We believe that this option, on average, should produce lower errors (w.r.t. $\Gamma^\uparrow$) and much more stable behaviors when compared to the metrics \cite{merge_farlocca, merge_farlocca_2, merge_wass}; 2) if there are multiple maps between metric merge trees giving the minimal interleaving value, clearly approximation errors occur only when all the cost-minimizing couplings force some unwanted behavior.

We leave to future works the assessment of the properties of the approximation scheme with $\mid f(father(x))-f(x)\mid u_x$.
We just point out that $\mid f(father(x))-f(x)\mid u_x$ in general in neither a lower or an upper bound to the cost of deleting $x$ and that is why we are not considering it in our theoretical investigation. Lastly, replacing $\mid g(r_G)-f(x)\mid u_x$ with the exact the deletion cost, would turn the linear optimization problem into a polynomial one, which we believe would make the computational cost too high.      

\subsection{Linearization}

At this point we have introduced a set of linear constraints, needed to optimize a non linear function (either $\Gamma^\uparrow$ or  $\Gamma^\downarrow$) of the form $\min_{\mathcal{V}\in\mathcal{K}}\max_i F_i(\mathcal{V})$ for some real-valued functions $F_i$ which are linear in $\mathcal{V}$. 
We can turn this into a linear optimization problem by exploiting a standard trick, introducing auxiliary variables and with additional constraints. 

Suppose we need to find $\min_{s} \max (f(s),g(s))$, with $f,g$ real valued functions; we then introduce the variable $u$, with the constraints $u\geq f(s)$ and $u\geq g(s)$ and solve the problem $\min_{s,u\geq f(s),u\geq g(s)} u$. 
We want to use this procedure to compute $\min_{\mathcal{V}\in \mathcal{K}} \Gamma^{\uparrow}(\mathcal{V})$ and $\min_{\mathcal{V}\in \mathcal{K}} \Gamma^{\downarrow}(\mathcal{V})$. We write down the details only for $\min_{\mathcal{V}\in \mathcal{K}} \Gamma^{\uparrow}(\mathcal{V})$, the case $\Gamma^{\downarrow}(\mathcal{V})$ follows easily.

Given $x\in E_T$ we define $F_x^1=\sum_{y} a_{x,y}\parallel C_{x,y}\parallel_\infty$, $F_x^2=(g(r_G)-f(x))u_x$, and $F_x^3=A_x$. Given $y\in E_G$, instead, we set $F_y^1=(f(r_T)-g(y))u_y$, and $F_y^2=A_y$. Analogously, we have $F_{r_T}^1=\mid f(r_T)-g(r_G) \mid $.
Having fixed a total ordering on $V_T=\{a_0,\ldots, a_n\}$ and $V_G=\{b_0,\ldots, b_m\}$, respectively, we call $\mathcal{F}$ the vector containing all the components $F_x^i$ and $F_y^j$ of all the points taken in the chosen order: 
\begin{equation}
\mathcal{F}:=\left( F_{a_1}^1,F_{a_1}^2,F_{a_1}^3,\ldots,F_{a_i}^1,F_{a_i}^2,F_{a_i}^3,\ldots,F_{b_1}^1,F_{b_1}^2,\ldots,F_{b_m}^1,F_{b_m}^2 \right)=
\left( F_i\right)
\end{equation}

Similarly, for every $x$, we order the elements of the set $B_x$, so that $B_x=\left(\ldots, B_x^h,\ldots \right)$; then we set

\begin{equation}
\mathcal{B}:=\left( B_{a_0}^1,\ldots,B_{a_0}^{h_0},\ldots,B_{a_i}^1,\ldots,B_{a_i}^{h_i},\ldots,B_{b_0}^1,\ldots,B_{b_0}^{t_0},\ldots,B_{b_j}^1,\ldots,B_{b_j}^{t_j}\right)=\left( B_i \right)
\end{equation}
 
So we introduce the real valued auxiliary variable  $z$ and add the following constraints to the ones presented in \Cref{sec:constr_interl}:
\begin{itemize}
\item[(5)] $z\geq F_i$ for all $i$;
\item[(6)] $z\geq B_i$ for every $i$;
\end{itemize}
and then solve $\min_{z,\mathcal{V}} z$ with all these constraints.
 We stress again that $F_i$ and $B_i$ are linear in $\mathcal{V}$; so the final problem is linear with mixed binary valued variables (actually all but one variable are binary). In case of $\Gamma^{\downarrow}$ we repeat the same operations, omitting the constraints in (6).

\subsection{Bottom-Up Algorithm and Computational Complexity}
\label{sec:bottom-up}

In this section the results obtained in  \Cref{sec:decomposition_interl} and the formulation established in the previous parts of \Cref{sec:LIP} are used to obtain the algorithm implemented to approximate the metric $d_I$ between  merge trees.
We need to introduce some last pieces of notation in order to describe the \virgolette{bottom-up} nature of the procedure.

Given $x\in V_T $, define $\len(x)$ to be the cardinality of $\{v \in V_T\mid x\leq v\leq r_T\}$ and $\len(T)= \max_{v\in V_T}\len(v)$; similarly $\lvl(x)= \len(T)-\len(x)$ and $\lvl_T(n)=\{v\in V_T\mid \lvl(v)=n\}$

The key property is that: $\lvl(x)> \lvl(v)$ for any $v\in sub_T(x)$.
Thus, for instance, if $W^\uparrow_{x,y}$ is known for any $x\in \lvl_T(n)$ and $y\in \lvl_{G}(m)$, then  for any $v$, $w$ in $V_T$, $V_{G}$ such that $\lvl(v)<n$ and $\lvl(w)<m$, $W^\uparrow_{v,w}$ is known as well. 
We write down Algorithm \ref{alg:bottomup_interl}, which refers to the computation of $\min \Gamma^{\uparrow}$. Note that thanks to constraints (4) this recursive procedure is always guaranteed to provide a special coupling.

In case of $\Gamma^{\downarrow}$ we repeat the same algorithm, omitting the constraints in (4) and (6) and with $W^\downarrow_{x,y} =\min_{z,\mathcal{V}}z$ being a lower bound for $d_I(T_x,G_y)$.

\begin{algorithm}[H]
\SetAlgoLined
\KwResult{Upper bound for $d_I(T,G)$ }
 initialization: $N=\len(T)$, $M = \len(G)$, $n=m=0$\;
 \While{$n\leq N$ or $m\leq M$}{
	 \For{$(x,y)\in V_T\times V_{T'}$ such that $\lvl(x)\leq n$ and $\lvl(y)\leq m$}{
	 	 Calculate $H_{x,y}$;\\
	 	 Calculate $W^\uparrow_{x,y} =\min_{z,\mathcal{V}}z$ subject to constraints (1)-(6), using $W^\uparrow_{x',y'}$ as upper bound for $\parallel C^o_{x',y'}\parallel_\infty$, for all couples $(x',y')$ already considered\;
	 } 
	$n=n+1$;
	$m=m+1$;
   
 }
\Return $\min_{(x,y)\in V_T\times V_G} \max \{H_{x,y},W^\uparrow_{x,y}\}$
 \caption{Bottom-Up Algorithm.}
 \label{alg:bottomup_interl}
\end{algorithm}

\begin{prop}[Computational Cost]\label{prop:complexity}
Let $T$ and  $T'$ be two merge trees with full binary tree structures with $\dim(T)=\#E_T = N$ and $\dim(T')=M$.

Then:
\begin{itemize}
\item to compute $W^\uparrow_{x,y}$ for every $(x,y)\in E_T\times E_{T'}$ with Algorithm \ref{alg:bottomup_interl}  we need to solve $O(M\cdot N)$ mixed binary linear optimization problems, with $O(M\cdot N)$ variables and $O(M\cdot \log_2(M)+ N\cdot \log_2(N))$ linear constraints;
\item to compute $W^\downarrow_{x,y}$ for every $(x,y)\in E_T\times E_{T'}$ with Algorithm \ref{alg:bottomup_interl}  we need to solve $O(M\cdot N)$ mixed binary linear optimization problems, with $O(M\cdot  N)$ variables and $O(\log_2(M)+ \log_2(N))$ linear constraints.
\end{itemize}

\begin{proof}
In a full binary tree structure, at each level $l$ we have $2^l$ vertices. Let $L=\len(T)$ and $L'=\len(T')$. We have that, for any vertex $v\in V_T$ at level $l$, the cardinality of the path from $v$ to any of the leaves in $sub_T(v)$ is ${L-l}$ and the number of leaves in $sub_T(v)$ is $2^{L-l}$. 
The number of vertices in $sub_T(v)$ is instead $2^{L-l+1}-1$ and the number of edges $2^{L-l+1}-2$.

Consider $v\in V_T$ at level $l$ and $w\in V_{T'}$ at level $l'$.  According to \Cref{sec:constr_interl} to compute $\Gamma^\uparrow$ with $T_v$ and $T'_w$  we need $(2^{L-l+1}-2)\cdot (2^{L'-l'+1}-2) + 2^{L-l+1}-2+ 2^{L'-l'+1}-2$ binary variables and $ 2^{L-l}+2^{L'-l'}+ 2\cdot (2^{L-l+1}-2)+ 2\cdot (2^{L'-l'+1}-2)$ linear constraints (if we consider constraints (1)-(3) in \Cref{sec:constr_interl}) in the nonlinear problem. Constraint (4) adds just $1$ further constraint, so we can ignore it. To linearize the problem we need to add 
$1$ real valued auxiliary variable and $3\cdot (2^{L-l+1}-2)+2\cdot (2^{L'-l'+1}-2)$ constraints to take into account for the set $\mathcal{F}$ and additional:
\[\sum_{k=0}^l 2^k\cdot (2^{L-l-k+1}-1)+\sum_{s=0}^{l'} 2^s\cdot (2^{L'-l'-s+1}-1)
\] 
 constraints for the set $\mathcal{B}$.
Simplifying we have: $(l+1)\cdot (2^{L-l+1}-1) + (l'+1)\cdot (2^{L'-l'+1}-1)\leq O(2^L + 2^{L'})$ constraints for $\mathcal{B}$.

Putting the things together, for the linearized problem we have:
\begin{equation}
(2^{L-l+1}-2)\cdot (2^{L'-l'+1}-2) + 2^{L-l+1}-2+ 2^{L'-l'+1}-2
\end{equation}
 binary variables. And:
\begin{equation}
(2^{L-l}+2^{L'-l'}+ 5\cdot (2^{L-l+1}-2)+ 4\cdot (2^{L'-l'+1}-2) + l\cdot (2^{L-l+1}-1) + l'\cdot (2^{L'-l'+1}-1)
\end{equation}
 linear constraints.

Thus, to minimize $\Gamma^\uparrow$, we need to solve $(2^{L+1}-1)\cdot (2^{L'+1}-1)$ linear binary optimization problems, each with  $2^L\cdot 2^{L'}+O(L+ L')$ variables and $O(2^L+ 2^{L'})$ constraints. 
Instead, for $\Gamma^\downarrow$, we need to solve $(2^{L+1}-1)\cdot (2^{L'+1}-1)$ linear binary optimization problems, each with  $2^L\cdot 2^{L'}+O(L+ L')$ variables and $O(L+ L')$ constraints.
Substituting $L = \log_2(N)$ and $L' = \log_2(M)$ in these equations gives the result.
\end{proof}
\end{prop}

\begin{rmk}[Computational Complexity]\label{rmk:comp_cost}
We make some comments about the computational cost of Algorithm \ref{alg:bottomup_interl} as described in 
\Cref{prop:complexity}. 
First we note that the classical edit distance between unlabelled and unordered trees obtained with BLP \citep{TED}, 
an be computed by solving $O(N\cdot M)$ BLP problems with $O(N\cdot  M)$ variables and $O(\log_2(M)+ \log_2(N))$ constraints - $O(N+M)$ if we count also the constraints restricting the integer variables to $\{0,1\}$, as the authors of \cite{TED} do. 

Thus the approximating procedures we developed are not far from the same computational complexity of the widely used tree edit distance: $\min W^\uparrow_{x,y}$ differs by some linear factors while $\min W^\downarrow_{x,y}$ has the same complexity. 

This leads us to the following consideration: the couplings we define in the present work share many similarities with the mappings defined for edit distances, in particular with the ones in \cite{pegoraro2024finitely}, which generalize the classical tree edit distance mappings. 
For the tree edit distance, there exists polynomial time approximation algorithms \citep{zhang1996constrained}, relying on mappings with additional constraints. In \cite{pegoraro2024finitely} it is shown that such constraints are also compatible with the \virgolette{generalized} mappings therein defined. Thus we strongly believe that the approach 
of \cite{zhang1996constrained} could be adapted to our coupling and could potentially lead to a polynomial time approximation of the interleaving distance. We point out that the most efficient metrics to compare merge trees \cite{merge_farlocca, merge_farlocca_2, merge_wass} are (unstable) edit distances which rely on the already mentioned polynomial approximation \citep{zhang1996constrained} of the classical edit distance or on other similar variations \citep{selkow1977tree, jiang1995alignment}. Lastly, it is hard to compare our numerical schemes with the one found in \cite{interl_approx}, as their algorithm is not written in terms of MBLP (and it is not an approximation). On top of that, to our knowledge, their solution has not been used in any data analysis scenario or simulation and a public implementation of their algorithm is not available. We just report that \cite{TED} finds massive performance improvements when resorting to BLP to compute the edit distance, compared to other (exact) computational solutions. Those advantages may be heavily dependent on the performances of solver employed. 
\end{rmk}

\section{Error Propagation}
\label{sec:errors}

We make a brief observation to take care of the interactions arising between the approximations of the interleaving distance defined in \Cref{sec:obj_fn} and the bottom-up procedure proposed in \Cref{sec:bottom-up}.

Consider the setting of \Cref{sec:obj_fn} and suppose that, instead of having computed the optimal couplings $C_{x,y}$, we have some approximations $C'_{x,y}$ - as it is the case for $W^\uparrow_{x,y}$ and $W^\downarrow_{x,y}$ - with an error bound $e_{x,y}$ such that $\mid \parallel C_{x,y} \parallel_\infty - \parallel C'_{x,y} \parallel_\infty \mid < e_{x,y}$.
We immediately see that the errors $e_{x,y}$ affect only the components of the form $ \sum_{y} a_{x,y}\parallel C_{x,y}\parallel _\infty$ as these are the only parts of the optimization procedure which depend on $\parallel C'_{x,y} \parallel_\infty$.
Moreover, the potential errors occurring in $\Gamma^\uparrow$ and $\Gamma^\downarrow$, appear at the level of $B_x$ and $u_x$ (see \Cref{sec:approx}). But  $u_x$ and $B_x$ do not depend on 
$\parallel C'_{x,y} \parallel_\infty$ or on $u_{x'}$ and $B_{x'}$ for some other $x'$. Which implies that 
the only interaction between errors is of the form $\max \{a,b\}$, but they never aggregate in any way at any time in the objective functions.  Which means that errors do not propagate exploding in size: the final error in the algorithm presented in \Cref{sec:bottom-up} is the maximum of all the \virgolette{independent} errors occurring at every iteration.

\section{Simulation Study}
\label{sec:simulation}

In this section we test our approximations versus another method to approximate $d_I$ recently proposed by \cite{curry2021decorated}, relying on the work of \cite{merge_intrins}.
The approximation proposed by \cite{curry2021decorated} turns the unlabeled problem of the interleaving distance between merge trees into a labeled interleaving problem by proposing a suitable set of labels. The optimal labeling would give the exact value of the interleaving distance, but, in general, this procedure just returns an upper bound. We call $d_{lab}$ the approximation obtained with the labeled method proposed by \cite{curry2021decorated} and $d_u$ and $d_l$ respectively our approximation from above (i.e. with $\Gamma^\uparrow$) and below (i.e. with $\Gamma^\downarrow$).

For any fixed $i$, we generate a couple of point clouds $C^k_i=\{(x^k_j,y^k_j)\mid j=1,\ldots, n_i\}\subset\mathbb{R}^2$, with $k=1,2$, according to the following process:

\begin{align*}
x^k_j \sim^{iid} &\mathcal{N}(5,\sigma_{x,k})\text{ }j=1,\ldots,n_i \\
y^k_j \sim^{iid} &\mathcal{N}(5,\sigma_{y,k})\text{ }j=1,\ldots,n_i \\
\sigma_{x,k} \sim &\mathcal{N}(3,1)\\
\sigma_{y,k} \sim &\mathcal{N}(3,1).
\end{align*}

Note that $n_i$ regulates the number of leaves in the trees (which we fix before sampling $C^k_i$).
From $C^k_i$ we obtain the single linkage hierarchical clustering dendrogram $T_{C^k_i}$ (that is, the merge tree representing the Vietoris Rips filtration of $C^k_i$). And then compute $d_u(T_{C_i^1},T_{C_i^2})$, $d_l(T_{C_i^1},T_{C_i^2})$
and $d_{lab}(T_{C_i^1},T_{C_i^2})$. The distance $d_{lab}$ is computed with the code available at \url{https://github.com/trneedham/Decorated-Merge-Trees}, while $d_u$ and $d_l$ are computed via the procedure described in \Cref{sec:bottom-up}. 
For each $n_i\in\{2,\ldots,15\}$ we sample $100$ pairs of point clouds and then compute the relative error $(d_{lab}-d_u)/d_u$ and $(d_{l}-d_u)/d_u$. 

We repeat the same experiment, this time with 
$n_i\in\{100,101,102\}$. Since in this case $d_u$ requires too much time to be computed exactly, we exploit \Cref{prop:pruning} and consider the smallest $\varepsilon>0$ such that $P_\varepsilon(T_{C^1_i})$
and $P_\varepsilon(T_{C^2_i})$ have fewer or equal then $15$ leaves.
We then call $d_{opt}(T_{C_1},T_{C_2})=\max\{\varepsilon/2, d_u(P_\varepsilon(T_{C_1}),P_\varepsilon(T_{C_2}))\}$.

In general we have the following inequality:
\[
d_I(T,G)\leq d_{opt}(T,G),d_{lab}(T,G). 
\]

\begin{figure}
    \begin{subfigure}[c]{0.49\textwidth}
    	\centering
    	\includegraphics[width = \textwidth]{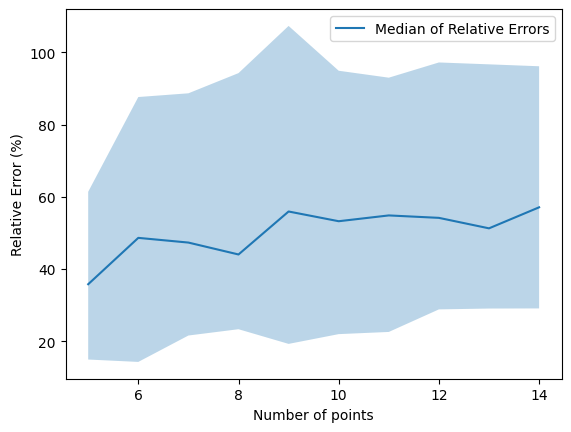}
		\caption{Median of the error percentage as a function of the number of leaves, with the shaded region being the central quartiles.}
		\label{fig:diff}
	\end{subfigure}
    \begin{subfigure}[c]{0.49\textwidth}
    	\centering
    	\includegraphics[width = \textwidth]{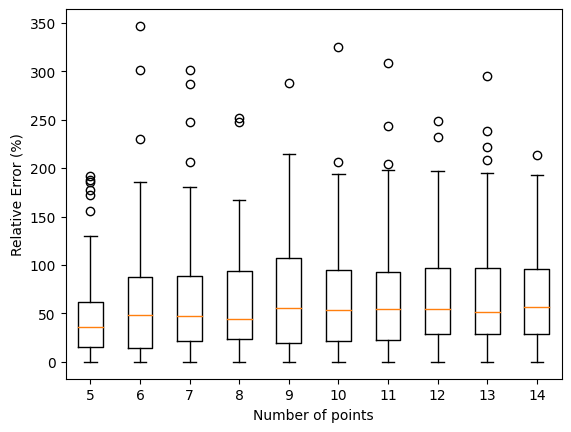}
		\caption{Boxplot of the error percentage as a function of the number of leaves.}
		\label{fig:mean}
	\end{subfigure}

    \begin{subfigure}[c]{0.49\textwidth}
    	\centering
    	\includegraphics[width = \textwidth]{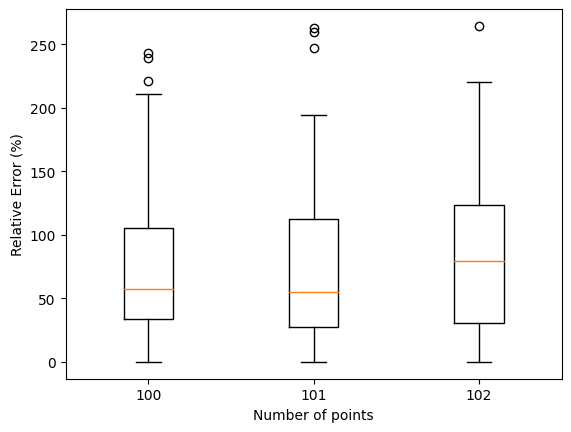}
		\caption{Boxplot of the error percentage as a function of the number of leaves, with the number of leaves being big (and so resorting to $d_{opt}$).}
		\label{fig:big_n}
	\end{subfigure}	
    \begin{subfigure}[c]{0.49\textwidth}
    	\centering
    	\includegraphics[width = \textwidth]{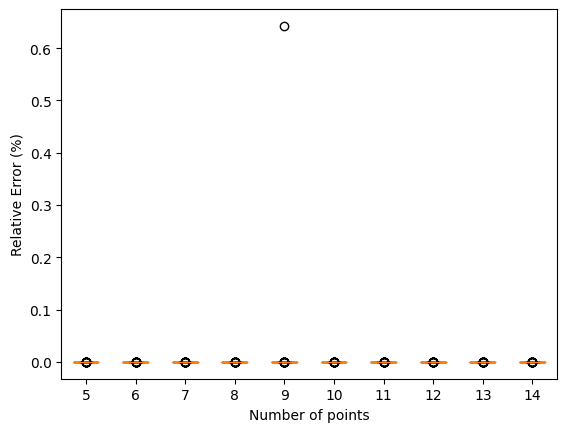}
		\caption{Boxplot of the differences between the upper and lower bounds to the interleaving distance, as a function of the number of leaves.}
		\label{fig:max}
	\end{subfigure}	

\caption{Descriptive statistics of the error percentages obtained with the approximations $d_l$, $d_u$, $d_{opt}$ and $d_{lab}$, as a function of the number of leaves, with $n_i\in\{2,\ldots,15\}$ - \Cref{fig:diff}, \Cref{fig:mean}, \Cref{fig:max} - and $n_i\in\{100,01,102\}$ - \Cref{fig:big_n}.}
\label{fig:simulation_1}
\end{figure}

The results of the simulations can be seen in \Cref{fig:simulation_1}. 
Looking at \Cref{fig:diff}, \Cref{fig:mean}, \Cref{fig:big_n} we see that, in the context of our data-generating process, the approximation given by $d_{lab}$ is very unreliable, producing a median error which around $50-60\%$, but with some outliers which are completely off from the values obtained with $d_u$, with errors of more than $3$ times the actual value.
\Cref{fig:max} instead shows that, with our data-generating pipeline, the gap between $d_u$ and $d_l$ is almost always $0$, which means that we are exactly computing the interleaving distance. The computations are carried out on a 2016 laptop with Intel(R) processor Core(TM) i7-6700HQ CPU @ 2.60GHz, 4 cores (8 logical) and 16 GB of RAM.
The employed MBLP solver is the freely available IBM CPLEX Optimization Studio 12.9.0.


As a conclusive remark, we say that the computational advantages of the labeled approach are immense and potentially adequate even for real time applications, but from our simulation we see that the results need to be taken with care, for the approximation produced is not always good. 
On the other hand, with the present implementation, the computational cost of our approach becomes prohibitive quite quickly as the number of leaves in the trees increases, even if there might be situations - like the one in this simulation - in which the approximation scheme we used for 
$n_i\geq 100$ can produce good estimates.
We think that interactions between the two approaches could lead to substantial speed-ups: being able to fix the value of some variables in our algorithm exploiting the labeling scheme by \cite{curry2021decorated} could greatly reduce the dimensionality of the problem and thus its computational cost.

\section{Discussion}
\label{sec:discussion_interl}

In this manuscript we propose a graph-matching approach to the interleaving distance between merge trees, trying to understand this metric by means of matchings between unlabeled combinatorial objects rather than with $\varepsilon$-good maps between metric trees. 
The relationship between this kind of matchings, which we call couplings, and continuous maps between merge trees helps us in producing another formulation of the interleaving distance.
In the second part of the manuscript we exploit this formulation to obtain some novel properties of the interleaving distance, along with two approximation schemes relying on a dynamical linear binary programming approach.
We test this schemes in a simulation study against another approximation procedure recently published, producing for the first time some error estimates for this second procedure. The algorithms we produce have some drawbacks in terms of computational cost, which limits their usage to small-data scenarios, but still show that the reliability of the other implementation available to work with this distance is often not very good. 

There are many directions which we would like to investigate in the future. We would like to further test the lower and upper bounds we propose in relationship to the third approximation scheme presented in \Cref{sec:approx}, perhaps forcing situations in which the lower and upper bound are far apart. Further testing would also benefit the understanding of upper bound we obtain via pruning. Similarly, we would like to compare the approximation of the interleaving distance with other distances proposed for merge trees, to assess the different \virgolette{unstable} behaviours.
On top of that, as we claimed in the end of \Cref{sec:map_couplings}, we believe that the relationship between a coupling $C$ and the map $\alpha^C$ can be used to obtain yet another definition of the interleaving distance using only continuous monotone maps between metric merge trees.

Lastly, the global modeling of the matching between the trees offered by couplings, opens up the door also to optimal couplings which also satisfy some additional local optimality criterion. This could be used to obtain some useful parametrizations of particular optimal maps between merge trees, leading to a better understanding of the geodesic structure of the space of merge trees with the interleaving distance.

\section*{Acknowledgments}
This work was carried out as part of my PhD Thesis, under the supervision of Professor Piercesare Secchi.

\section*{Declaration of competing interest}
The author declare that he has no known competing financial interests or personal relationships that could have
appeared to influence the work reported in this paper.

\section*{Data availability}
No real-world or benchmark data was used for the research described in the article. The generative process for the simulated data sets is described in the manuscript.

\section*{Code}

The code is available upon request. 
At the present stage it lacks proper documentation and commenting.



\appendix

\section{Proofs}\label{sec:proofs}

\bigskip\noindent
\underline{\textit{Proof of}  
\Cref{prop:alpha_order}.}

\smallskip\noindent

\begin{itemize}
\item The thesis is obvious if $x,x'\in\pi_T(C)$;
\item if $x\in \pi_T(C)$ and $x'\in D$, then $w\geq \chi(x')$, with $(x,w)\in C$, otherwise (C3) is violated since $\#\Lambda(x)>1$, and so $\alpha_C(x)\geq \alpha_C(x')$;
\item if $x'\in \pi_T(C)$ and $x\in D$, then $\chi(x)\leq w$ with $(x',w)\in C$;
\item lastly, suppose $x,x'\in D$ and consider the following cases:
	\begin{itemize}
	\item if $\#\Lambda(x)=\#\Lambda(x')=0$ then we have $\varphi(x)= \varphi(x')$ and thus $f(x)+\frac{1}{2}(f(\varphi(x))-f(x))\leq f(x')+\frac{1}{2}(f(\varphi(x'))-f(x'))$. This entails
 $g(\alpha_C(x))\geq g(\alpha_C(x'))$;
 	\item if $\#\Lambda(x),\#\Lambda(x')>1$ then we have $\chi(x)\leq \chi(x')$ and the thesis clearly follows;
	\item suppose lastly that $\#\Lambda(x)=0$ and $\#\Lambda(x')=1$; then $\eta(x)\leq \chi(x')$.	
	\end{itemize}
\end{itemize}

\hfill$ \blacksquare $

\bigskip\noindent
\underline{\textit{Proof of}  
\Cref{lem:sequences}.}

\smallskip\noindent

Given $x\in \mathbf{T}$, if $\max\pi_T(C)\in \{v\in \mathbf{T} \mid v\geq x\}$, we have a well defined upper extreme $u_C(x)$. Moreover, since $\{v\in \mathbf{T} \mid v\geq x\}$ is totally ordered, $u_C(x)$ is  unique. If $\max\pi_T(C)\notin \{v\in \mathbf{T} \mid v\geq x\}$, i.e. $\LCA(x,\max\pi_T(C))>\max\pi_T(C)$, then $u_C(x)=+\infty$.

Now consider $x\in L_T$. Since $sub_T(x)=\{x\}$ then $x\notin U$. Thus also $l_C(x)$ is well defined.
On top of that, since by hypotheses we are considering only vertices such that $x\notin \pi_T(C)\bigcup D$, $x\neq l_C(x)\neq u_C(x)$. Thus $[l_C(x),u_C(x)]$ is a non degenerate sequence of edges.

Suppose $\{v\in U \mid l_C(x)\leq v \leq x \} = \emptyset$
and $v',v''\in \max \{v \leq x \mid v \in \pi_T(C)\cup D\}$. Clearly $[v',x]\bigcap [v'',x] \bigcap V_T \neq \emptyset $; so consider $p \in [v',x]\bigcap [v'',x] \bigcap V_T$. We know $p\notin U$, thus $p \in \pi_T(C)\cup D$, which is absurd since $v',v''<p$.

Lastly, suppose $x\in U$ (and so $\#\Lambda(x)=1$) and $\max \{v\leq x \mid v \in \pi_T{C}\cup D\}\subset D$. Then  $\#\Lambda(l_C(x))\neq 1$ for any $l_C(x)$. 
Let $\{v\}=\Lambda(x)$ and consider $[v,x]$. If $v'\in U$ for all $v<v'<x$, we are done. Clearly there cannot be vertices $v'$ with $v<v'<x$ which are in $\pi_T(C)$. So suppose there is $v'\in D$, with  $v<v'<x$.
Since $v\in \Lambda(v')$, we have $\#\Lambda(v')>0$, but then $\#\Lambda(v')>1$ which means $v,v''\in \Lambda(v')$ for some $v''<v'$. Clearly in $[v'',x]$ there can be no vertex 
contained in $\pi_T(C)$ apart from $v''$. Thus $\#\Lambda(x)>1$, which is absurd. 

\hfill$ \blacksquare $

\bigskip\noindent
\underline{\textit{Proof of}  
\Cref{prop:extension}.}

\smallskip\noindent

We need to check continuity.
Consider $x_n\rightarrow x$ in $\mathbf{T}$. We know that their order relationships is preserved by $\alpha_C$ and $\alpha_C^\varepsilon$. On top of that $f(x_n)\rightarrow f(x)$ yields $f(x_n)+\varepsilon\rightarrow f(x)+\varepsilon$ and the result follows.

\hfill$ \blacksquare $

\bigskip\noindent
\underline{\textit{Proof of}  
\Cref{lem:diseq_interl}.}

\smallskip\noindent

We know that $\Lambda(x),\Lambda(y)>1$. Thus $x,y$ are either deleted or coupled. Note that $\chi(x)\geq y$ and $\chi(y)\geq x$.
If both of them are coupled then $(x,y)\in C$. Suppose $x$ is coupled - with $\delta(y)$ - and $y$ is deleted. Then $cost(y)=\mid f(\chi(y))-g(y)\mid < \varepsilon$
 and 
$cost(x)=\mid f(x)-g(\delta(y))\mid < \varepsilon$. Since $\delta(y)\geq y$ and $\chi(y)\geq x$ we obtain the thesis.

\hfill$ \blacksquare $

\bigskip\noindent
\underline{\textit{Proof of}  
\Cref{teo:eps_to_C}.}

\smallskip\noindent

\begin{itemize}
\item $\alpha_C^\varepsilon$ is continuous by \Cref{prop:extension}.

\item (P1) holds by \Cref{prop:extension}. 

\item Now we prove (P2). 
Suppose we have $\alpha_C^\varepsilon(v)<\alpha_C^\varepsilon(v')$. And let $x=\LCA(v,v')$.

We can suppose $\#\Lambda(x)>1$ and $\LCA(\Lambda(x))=x$; otherwise at least one between $\varphi(v)\geq x$ and $\varphi(v')\geq x$ holds. Suppose the second one holds, then $\#\Lambda(v')=0$ and $\varepsilon\geq (f(\varphi(v'))-f(v'))/2 $ and $\varphi(v')\geq x$. Thus, (P2) holds. The same if the first one holds.

Now we show that if $\#\Lambda(x)>1$ we can find $a,b\in V_T$ such that:
\begin{itemize}
\item $x=\LCA(a,b)$;
\item $(a,a'),(b,b')\in C$;
\item $\alpha_C(v)\geq \alpha_C(a)$ and $\alpha_C(v')\geq \alpha_C(b)$.
\end{itemize}

Note that, in this case, upon calling $y= \LCA(a',b')$ we have: $\mid f(x)-g(y) \mid \leq  \varepsilon$ by \Cref{lem:diseq_interl}, $\alpha^\varepsilon_C(v')\geq \{\alpha_C(v),\alpha_C(a), \alpha_C(b)\}$ and so  
$\alpha_C(v')\geq y$. Which means that $f(v')+\varepsilon \geq g(y)$. Thus $ f(x)-f(v') \leq 2\varepsilon $.

We enumerate all the possible situations for $v$ - clearly the same hold for $v'$:
\begin{itemize}
\item $v\in \pi_T(C)$: then $a=v$;
\item $\Lambda(v)=0$ then $a=v''$ with $(v'',\eta(v))\in C$; by hypothesis, $\varphi(v)\leq x$ and $v''<x$; 
\item $a\in\Lambda(v)$ if $v\notin \pi_T(C)$ and $\#\Lambda(v)>0$.
\end{itemize}

\item Now we prove (P3). If $w \notin Im(\alpha_C^\varepsilon)$ then $w \in D_C^G$. In fact if $(x,w')\in C$ with $w'<w$, then $[w',w]\subset Im(\alpha_C^\varepsilon)$, since there are $l(w)\geq w'$ and $u(w)\leq r_{G}$. But then we know that there is 
$w''=\min \{y\in V_G \mid y> w \text{ and }y\notin D_C^G\}$ with $g(w'')-g(w)\leq \varepsilon$.
\end{itemize}

\hfill$ \blacksquare $

\bigskip\noindent
\underline{\textit{Proof of}  
\Cref{lem:lemma_0}.}

\smallskip\noindent

We know that $w:=\alpha(v)\leq w':=\alpha(v')$, thus 
$\max\{y \in V_G \mid y\geq w\}\leq \max\{y \in V_G \mid y\geq w'\}$ and the thesis follows.

\hfill$ \blacksquare $

\bigskip\noindent
\underline{\textit{Proof of}  
\Cref{teo:C_to_eps}.}

\smallskip\noindent

We build $C$ by subsquently adding couples starting from an empty set. The proof is divided in sections which should help the reader in following the various steps.

\subsection{Leaves of $T$}\label{sec:proof_leaf_T_interl}
In this section we take care of the leaves of the merge tree $T$.

\subsubsection{Selecting the Coupled Leaves}
\label{sec:proof_L_T_interl}
We consider the following set of leaves:
\begin{equation}
\mathcal{L}_T=\{v\in L_T\mid \nexists v'\in L_T \text{ such that }\alpha(v)< \alpha(v')\}
\end{equation}
We give a name to the condition:
\begin{itemize}
\item[$(a)$] $\nexists v'\in L_T \text{ such that }\alpha(v)< \alpha(v')$
\end{itemize}
so that we can more easily use it during the proof.
Note that we can avoid treating the case $\alpha(v)=\alpha(v')$ thanks to $(G)$.

The set $\mathcal{L}_T$ is the set of leaves which will coupled by $C$, while all other leaves will be deleted: we add to $C$ all the couples of the form $(v,\phi(v))$ with $v\in \mathcal{L}_T$. We characterize those couples with the following proposition.

\begin{lem}\label{prop:L_T_interl}
Given $v,v'\in \mathcal{L}_T$, then $\phi(v)\geq \phi(v')$ if and only if $v=v'$. Moreover, for every $v'\in L_T$ such that $(a)$ does not hold, there is $v\in\mathcal{L}_T$ such that $\alpha(v)<\alpha(v')$.

\begin{proof}
The first part of the proof reduces to observing that $\phi(v)\leq \phi(v')$ if and only if $\alpha(v)\leq \alpha(v')$.

Now consider $v'\in L_T$ such that $(a)$ does not hold.
We know there is $v_0$ such that $\alpha(v_0)<\alpha(v')$. If $v_0\in\mathcal{L}_T$ we are done, otherwise there is $v_1$ such that $\alpha(v_1)<\alpha(v_0)<\alpha(v')$. Note that $f(v_1)<f(v_0)$. Thus we can carry on this procedure until we find $v_i\in \mathcal{L}_T$. Note that $\arg\min_{V_T} f \in \mathcal{L}_T$, thus, in a finite number of step we are done. 
\end{proof}
\end{lem}

\subsubsection{Cost Bound on Couples}\label{sec:proof_cost_couples_interl}
Now we want to prove the following proposition which gives an upper bound for the cost of the couples added to $C$.

\begin{lem}\label{prop:cost_leaves_interl}
Given $v\in \mathcal{L}_T$, then $\mid f(v)-g(\phi(v))\mid \leq \varepsilon$.

\begin{proof}
Suppose the thesis does not hold. 
Since $g(\phi(v))\leq f(v)+\varepsilon$, contradicting the thesis means that we have $v\in\mathcal{L}_T$ such that: 
\begin{itemize}
\item[$(b)$] $g(\phi(v))+\varepsilon < f(v)$.
\end{itemize}
Let $w=\phi(v)$.
If $(b)$ holds, then $g(father(w))-g(w)>g(\alpha(v))-g(w)> 2\varepsilon$.
Let $v'=\psi(w)\leq \beta(w)$. Note that $f(v')<f(v)$. We have $\phi(v')\leq \alpha(v')\leq  \alpha(\beta(w))=s^{2\varepsilon}_G(w)$. But since $g(father(w))-g(w)>2\varepsilon$, we also have
$\alpha(v')\leq \alpha(v)$ with $v'\neq v$ which is absurd by \Cref{prop:L_T_interl}.
\end{proof}
\end{lem}

\subsubsection{Cost Bound on Deletions}\label{sec:proof_cost_deletions_interl}
In this step we prove the following proposition which gives an analogous bound to the one of \Cref{prop:cost_leaves_interl}, but for the deleted leaves of $T$.

\begin{lem}\label{prop:cost_deletions_interl}
Given $v\in L_T-\mathcal{L}_T$, then there exists $x>v$ such that:
\begin{itemize}
\item there is $v'<x$ such that $v'\in \mathcal{L}_T$;
\item $f(x)\leq f(v)+2\varepsilon$.
\end{itemize}

\begin{proof}
Since $(a)$ does not hold for $v$, we use \Cref{prop:L_T_interl} to obtain $v'\in \mathcal{L}_T$ such that $\alpha(v')< \alpha(v)$. But being $\alpha$ an 
$\varepsilon$-good map, we have $s_T^{2\varepsilon}(v')\leq s_T^{2\varepsilon}(v)$ which implies $f(\LCA(v,v'))\leq f(v)+2\varepsilon$. Thus $x=\LCA(v,v')$ ends the proof.
\end{proof}
\end{lem}

\Cref{prop:cost_deletions_interl} implies that, using the notation of the proposition, $\varphi(v)\leq x$. Then $f(v)<f(v')$ implies that $g(\phi(v))\leq f(v')+\varepsilon < f(v)+\varepsilon$. Thus $g(\eta(v))< f(v)+\varepsilon$. Since $f(x)\leq f(v)+2\varepsilon$ we have that the cost of deleting 
any $x'<x$ with $\#\Lambda(x')=0$ is less then $\varepsilon$.

\subsection{Leaves of $G$}\label{sec:proof_L_G_interl}
The result we need in this section is the following.

\begin{lem}\label{prop:G_leaves_interl}
Given $w\in L_G$, there exist $y\geq w$ such that:
\begin{itemize}
\item there is $w'=\phi(v)$ with $v\in\mathcal{L}_T$ and $w'<y$;
\item $g(y)\leq g(w)+2\varepsilon$.
\end{itemize}

\begin{proof}
Consider $\beta(w)$. 
Let $v\leq \beta(w)$ leaf. 
We have $\alpha(\beta(w))\geq \LCA(\alpha(v),w)$. If $v\in\mathcal{L}_T$ we are done for $\phi(v)\leq \alpha(v)$. If $(a)$ does not hold by \Cref{prop:L_T_interl} it means that there is $v'\in \mathcal{L}_T$ such that $\alpha(v')< \alpha(v)$. We are done since $g(\alpha(\beta(w)))= g(w)+2\varepsilon$.
\end{proof}
\end{lem}

Note that if $w<\phi(v)$ for some $v\in\mathcal{L}_T$, then, by \Cref{prop:G_leaves_interl}, $g(\phi(v))\leq g(w)+2\varepsilon$. In fact, using the notation of \Cref{prop:G_leaves_interl}, in this case we have $w'=y=\phi(v)$ by definition.
As in \Cref{sec:proof_cost_deletions_interl}, we have that the cost of the deletion of any $w$ such that $\#\Lambda(w)=0$ and $w\notin \pi_G(C)$ is at most $\varepsilon$.

\subsection{Internal Vertices}\label{sec:proof_internal_vert_interl}
Now we need to extend the coupling $C$ taking into account the internal vertices of $T$. We will do so after simplifying our merge trees in two different ways: first we remove all vertices which are deleted with $\#\Lambda(p)=0$ and then we take out all inessential internal vertices.

\subsubsection{Pruning}\label{sec:proof_pruning}
Let $T_0=T$ and $G_0=G$.
We define $T_1$ as the merge tree obtained from $T_0$ deleting the following set of vertices (and the corresponding edges): $x$ satisfying both $\#\Lambda(x)=0$ and $x\notin \pi_T(C)$.
Note that, for any $x\in V_T$ either there is $v\leq x$ with $v\in\mathcal{L}_T$ or for any leaf below $x$, we can apply \Cref{prop:cost_deletions_interl} and the consequential observations.

The tree $G_1$ is obtained from $G_0$ in an analogous way: 
any time we have $w\in V_G$ satisfying both $\#\Lambda(w)=0$ and $w\notin \pi_G(C)$, $w$ is deleted from $G_0$, along with the edge $(w,father(w))$.

Before proceeding we point out that, by construction, the leaves of $T_1$ are exactly $\mathcal{L}_T$.

\subsubsection{Restricting $\alpha$}
\label{sec:proof_alpha}
Thanks to \Cref{lem:lemma_0} we have that, anytime we delete some vertex in $G_0$ to obtain $G_1$ and that vertex is in the image of $\phi$, we are sure that also its counterimage is deleted from $T_0$.
Now, for every $v\in\mathbf{T}$, by construction we have that $\alpha(v)$ belongs to an edge removed from $G_0$ if and only if $\phi(v)$ is deleted from $G_0$. Let $v'=L(v)\in V_T$.  
Then $\phi(v')\leq \phi(v)$ and thus $\phi(v')$ is deleted as well. Which entails that $v'$ is deleted as well.
All of this, put together implies that we can restrict $\alpha$ to $\mathbf{T}_1$ (the metric tree obtained from $T_1$) and its image lies in $\mathbf{G}_1$ (obtained from $G_1$).

\subsubsection{$\varepsilon$-Good Restriction}
\label{sec:proof_alpha_eps}
We define $\alpha_1:=\alpha_{\mid \mathbf{T}_1}:\mathbf{T}_1\rightarrow \mathbf{G}_1$. 
We want to prove that $\alpha_1$ is still an $\varepsilon$-good map. Clearly $(P1)$ and $(P2)$ still hold upon restricting the domain. We just need to show $(P3)$. Suppose there is $w\in V_{G_0}$ such that $w\notin \alpha_1(\mathbf{T}_1)$. We distinguish between two cases: $(1)$ $w\notin \phi(\mathbf{T}_1)$ and $(2)$  
$w\in \phi(\mathbf{T}_1)$. Consider scenario $(1)$: 
$w\notin \phi(\mathbf{T}_1)$ clearly implies that $\#\Lambda(w)=0$ and so 
$w$ is deleted; scenario $(2)$ instead means that there is 
$\alpha(v)=\min \{\alpha(v')>w\mid v'\in\mathbf{T}_1\}$ with $L(\alpha(v))=\phi(v)= L(w)$. Clearly $\{w' \in \mathbf{G}_0\mid w'<w \}\bigcap \alpha(\mathbf{T}_1) =\emptyset$ and thus $v$ is a leaf of $T_1$, which means $v\in\mathcal{L}_T$. This also implies that $w\notin \alpha(\mathbf{T})$ and in particular $\alpha(v)=\min \{\alpha(v')>w\mid v'\in\mathbf{T}_0\}$, thus condition $(P3)$ is satisfied.

\subsubsection{Properties of $\phi:\mathbf{T}_1\rightarrow \mathbf{G}_1$ and Removal of Inessential Vertices} 
\label{sec:proof_order_2}
We know that $\phi:L_{T_1}\rightarrow L_{G_1}$ is injective. On top of that we have proved that, if $w\notin \alpha_1(\mathbf{T}_1)$ then $L(w)=\phi(v)$ for some $v\in\mathcal{L}_T$. Thus $\phi:L_{T_1}\rightarrow L_{G_1}$ is a bijection.

From now on we will ignore any vertex $v$ such that $\#child(v)=1$.
Formally, we introduce $T_2$ (and $G_2$) obtained from $T_1$ ($G_1$) removing all the vertices $v\in V_{T_1}$ such that $\#child(v)=1$ (similarly $w\in V_{G_1}$ such that $\#child(w)=1$): consider $\{v'\}=child(v)$. We remove $v$ from $V_{T_1}$ and replace the edges $(v',v)$ and $(v,father(v))$ with the edge $(v',father(v))$.
We do this operation recursively until no vertices such that $\#child(v)=1$ can be found.

\subsubsection{Coupling and Deleting the internal vertices}
\label{sec:proof_coupl_int}

We start with the following lemma.

\begin{lem}\label{lem:chi}
For every $x\in V_{T_2}$, and $y\in V_{G_2}$, we have $\mid x-\chi(x)\mid \leq \varepsilon$ and 
$\mid y-\chi(y)\mid \leq \varepsilon$.
\begin{proof}
Let $x\in V_{T_2}-L_{T_2}$. Then $x=\LCA(v_1,\ldots,v_n)$ with $v_1,\ldots,v_n$ being the leaves of $sub_{T_2}(x)$. Then $\alpha(x)>\alpha(v_i)$ for every $i$ and thus $\alpha(x)\leq w=\chi(x)=\LCA(\alpha(v_1),\ldots, \alpha(v_n))$. Clearly $x\leq \beta(w)$ for analogous reasons. Thus $\mid x-w\mid\leq \varepsilon$. 
For $y\in V_{G_2}$ we can make an analogous proof.
\end{proof}
\end{lem}

Let $v_1,\ldots, v_n\in L_{T_2}$ and  
$\phi(v_1),\ldots, \phi(v_n)\in L_{G_2}$.
Let $x=\LCA(v_1,\ldots, v_n)$ and $y=\LCA(\phi(v_1),\ldots, \phi(v_n))$. We know that $\alpha(x)\geq y$ and $\beta(y)\geq x$ and so $\mid f(x)-g(y)\mid \leq \varepsilon $.

Thus we proceed as follows: we order all the internal vertices of $T_2$ according to their height values and we start from the higher one (the root of $T_2$ - note that this in general is not the root of $T_0$), coupling it with the other root (thus $(C1)$ is verified).
Consider a lower $x$ and let $v_1,\ldots, v_n$ be the leaves of $sub_{T_2}(x)$ (and thus $x=\LCA(v_1,\ldots, v_n)$). Let $w=\LCA(\phi(v_1),\ldots, \phi(v_n))$. If the leaves of $sub_{G_2}(w)$ are $\phi(v_1),\ldots, \phi(v_n)$ we add the couple $(x,w)$ otherwise we skip $x$ - which will be deleted, with $\#\Lambda(x)>1$. Note that $w=\chi(x)$. Moreover, by \Cref{lem:chi}, $\mid f(x)-g(w)\mid \leq \varepsilon$.
Going from the root downwards we repeat this procedure for every $x$.

We clearly have:
\begin{itemize}
\item that (C3) is satisfied;
\item $cost((v,w))\leq \varepsilon$;
\item if $v\in V_{T_2}$ is deleted then $\#\Lambda(v)>1$ and $cost(v)=\mid f(v)-g(\chi(v))\mid \leq \varepsilon$;
\item if $w\in V_{G_2}$ is deleted then $\#\Lambda(w)>1$ and $cost(w)=\mid g(w)-f(\chi(w))\mid \leq \varepsilon$ by \Cref{lem:chi}. 
\end{itemize}

\subsection{Coupling Properties and Costs}

First we verify that $C$ is a coupling:
\begin{itemize}
\item[$(C1)$] verified. See \Cref{sec:proof_coupl_int}; 
\item[$(C2)$] verified. See  \Cref{sec:proof_leaf_T_interl} and \Cref{sec:proof_coupl_int}: we carefully designed the coupling so that no vertex is coupled two times;
\item[$(C3)$] verified. It is explicitly proven in \Cref{sec:proof_coupl_int};
\item[$(C4)$] In \Cref{sec:proof_order_2} we remove all vertices such that $\#child(x)=1$ to obtain$T_2$ and $G_2$; all the couples in $C$ are either leaves of vertices in $T_2$ and $G_2$. So, since all leaves of $T_2$ and $G_2$ are coupled its impossible to have a vertex $p$ belonging to any tree such that $\#\Lambda(p)=1$. 
\end{itemize}

Lastly we verify its costs:
\begin{itemize}
\item if $(x,y)\in C$ then $\mid f(x)-g(y)\mid \leq \varepsilon$ as verified in \Cref{sec:proof_cost_couples_interl} for $x\in L_T$ and in \Cref{sec:proof_coupl_int} for $x$ being an internal vertex;
\item if $x \in V_T\cap D$, with $\#\Lambda(x)=0$, it is verified that $\mid f(x)-\phi(x)\mid \leq 2 \varepsilon$ in \Cref{sec:proof_cost_deletions_interl}; if instead $y\in V_G\cap D$ we verify $\mid g(y)-\phi(y)\mid \leq 2\varepsilon$  in \Cref{sec:proof_L_G_interl}; in both cases we have a vertex lower that $x$ $(y)$ which, by construction, it is coupled. And the cost of the couple is less then $\varepsilon$. Thus the deletion of $x$ $(y)$ costs less than or equal to $\varepsilon$;
\item if $p \in D$, with $\#\Lambda(p)>1$, then 
we verify in \Cref{sec:proof_coupl_int} that the cost of this deletion is less than or equal to $\varepsilon$.
\end{itemize}
This concludes the proof.

\hfill$ \blacksquare $

\bigskip\noindent
\underline{\textit{Proof of}  
\Cref{teo:deco}.}

\smallskip\noindent

First, we prove the second part of the statement.

Given $C\in \mathcal{C}(T,G)$ optimal coupling such that $\{(r,r')\}=\max C$, we define 
$\mathcal{M}(C):= \max (C-
\{(r,r')\})$.
Clearly $\mathcal{M}(C)\in \mathcal{C}^*(T,G)$.
For any $(x,y)\in \mathcal{M}(C)$, consider $C^o_{x,y}\in \mathcal{C}_R^o(T_x,G_y)$ such that $\parallel C^o_{x,y}\parallel_\infty \leq \parallel C_{\mid (x,y)}\parallel_\infty$ with $C_{\mid (x,y)}:=\{(v,w)\in C \mid (v,w)<(x,y)\}$.  By hypothesis we know that $C^o_{x,y}$ exists for every $x,y$. Note that $C_{\mid (x,y)}\in\mathcal{C}_R(T_x,G_y)$.

We want to prove that the extension $C':= \{(r,r')\}\bigcup_{(x,y)\in \mathcal{M}(C)} C^{o}_{x,y}$ satisfies $\parallel C'\parallel_\infty \leq \parallel C\parallel_\infty$.
The set $C'$ clearly is a coupling since for every $x,x'\in\pi_T(\mathcal{M}(C))$, $sub_T(x)$
and $sub_T(x')$ are disjoint. And the same for $y,y'\in\pi_G(\mathcal{M}(C))$.
We need to consider the costs of different kinds of vertices separately. In particular we indicate with $(a)$ whenever we have $v\in V_T$ such that $v\leq x$ for some $x\in \pi_T(\mathcal{M}(C))$. 
If this condition does not hold we say that $(b)$ holds for $v$. The same definitions apply also for vertices in $G$. For instance $(a)$ holds for all vertices in $\pi_T(\mathcal{M}(C))$ and $\pi_G(\mathcal{M}(C))$. Similarly, if $(b)$ holds for $v$, then it holds also for all $v'>v$.

\begin{itemize}
\item Consider $v\in V_T$ such that $(a)$ holds for some $(x,y)\in \mathcal{M}(C)$; we have that $\chi_{C'}(v)\leq y$, $\varphi_{C'}(v)\leq x$ and $\eta_{C'}(v)\leq y$, so $cost(v)$ depends only on $C^o_{(x,y)}$.
Thus $cost_{C'}(v)=cost_{C^o_{(x,y)}}(v)\leq \parallel C_{\mid(x,y)}\parallel$ by construction; 
\item suppose now $(b)$ holds for $v\in V_T$; 
in this case we have by construction 
$\chi_C(v)=\chi_{C'}(v)$, $\Lambda_{C'}(v)=\Lambda_{C}(v)$ and thus $\varphi_C(v)=\varphi_{C'}(v)$. 
Suppose now $\Lambda_{C'}(v)=0$: $(a)$ holds for $\eta_{C'}(v)$ by means of some $y\in\pi_G(\mathcal{C}(M))$ and thus $g(\eta_{C'}(v))= \min_{w'\leq y}g(w')\leq g(\eta_{C}(v))$ (thanks to the properties of couplings in $\mathcal{C}^o_R(T_x,G_y)$). Since $\varphi_C(v)=\varphi_{C'}(v)$ we then have $cost_{C'}(v)\leq cost_C(v)$.
Lastly consider $\#\Lambda_{C'}(v)>1$. In this case we have $\chi_{C}(v)=\chi_{C'}(v)=r'$ and thus 
$cost_{C'}(v)= cost_C(v)$. Note that in all these situations, $cost_{C'}(v)$ does not depend on the chosen extension.
\end{itemize}

Exchanging the role of $T$ and $G$ we obtain the same results for the vertices of $G$.

The first part of the statement follows with analogous reasoning, by dropping the special coupling constraints and focusing on the vertices $v \in \pi_T(e(C^*))\text{ or } \#\Lambda_{e(C^*)}(v)>0$. In particular, the special couplings constraints are necessary only if $\#\Lambda_{e(C^*)}(v)=0$.

\hfill$ \blacksquare $

\bigskip\noindent
\underline{\textit{Proof of}  
\Cref{lem:pruning}.}

\smallskip\noindent

\begin{enumerate}
\item This is  due to the fact that all vertices of $P\varepsilon(T)$ are coupled, $P\varepsilon(T)$ is still a rooted tree and $(P)$ does not alter the order relationships of the remaining vertices;

\item since $(P)$ removes at most one leaf from the previous tree and adds no new leaves, it is clear that $L_{P_\varepsilon(T)}\subset L_T$. 
Consider now $l\in L_T$, $v_- = \max \{v'\in V_T \mid v\geq l \text{ and } f(v')< f(l)+\varepsilon\}$ and $v_+ = \min \{v'\in V_T \mid v\geq l \text{ and } f(v')\geq f(l)+\varepsilon\}$. Note that $(v_+,v_-)\in E_T$.
 By construction $v'=\arg\min_{v''\leq v_+} f(v'')$ is either the last point below $v_+$ to be deleted or is in $L_{P_\varepsilon(T)}$. Suppose it is the last point below $x_+$ to be removed by $(P)$.
By construction, when $(P)$ removes $v'$, $f(father(v'))-f(v')<\varepsilon$, but, since all other vertices below $x_+$ have been deletex $father(v')\geq x_-$ and so  
$f(father(v'))-f(v')\geq \varepsilon$, which is absurd. 
Exploiting this fact, also this point is proven;

\item consider $v\in V_T- V_{P_\varepsilon(T)}$; let $\Lambda_{C_\varepsilon}(v)=\{a_1,\ldots,a_n\}$. By construction  $v\geq x=\LCA(\Lambda_{C_\varepsilon}(v))$. We know that a vertex is removed from $T$ if, at a certain point along the recursive application of $(P)$, it becomes an order $2$ vertex or a small-weight leaf.
 Thus the vertex $x$ is not removed from $T$ unless $n\leq 1$. Which is absurd because then $\Lambda_{C_\varepsilon}(v)=\{x\}$. Thus  $\#\Lambda(v)\leq 1$.  As a consequence,  $v\in D$ if and only if $\#\Lambda(v)=0$. 
 
The vertex $\varphi_{C_\varepsilon}(v)$
is the first vertex $x$ above $v$ such that there is a leaf $v'$, $v'<x$, with 
$f(x)-f(v')\geq \varepsilon$. 
If $f(x)-f(v)\geq \varepsilon$ then by point $(2)$ there exist $v'\in L_{P_\varepsilon(T)}$ with $\LCA(v',v)<\varphi_{C_\varepsilon}(v)$ which is absurd;

\item consider two leaves $l,l'\in L_T$ with  $f(l)<f(l')$, with $l\notin V_{P_\varepsilon(T)}$ and $l'\in V_{P_\varepsilon(T)}$. On top of that suppose $l'<\varphi_{C_\varepsilon}(l)$.By point $(3)$ we have $f(\varphi_{C_\varepsilon}(l))-f(l')< \varepsilon$ which means that $l'$ is a small-weight leaf. Which is absurd;

\item combining points $(3)$ and $(4)$ we see that we only have deletions with $\#\Lambda_{C_\varepsilon}(v)=0$ and $2\cdot  cost_{C_\varepsilon}(c)= f(\varphi_{C_\varepsilon}(v))-f(v)\leq \varepsilon$.
\end{enumerate}

\hfill$ \blacksquare $

\bigskip\noindent
\underline{\textit{Proof of}  
\Cref{prop:pruning}.}

\smallskip\noindent

Consider $C\in\mathcal{C}(P_\varepsilon(T),P_\varepsilon(G))$. Then $C$ can be seen also as a coupling $C\in\mathcal{C}(T,G)$. Let $C_\varepsilon$ be as in \Cref{lem:pruning}.

We partition the vertices $V_T$ into three sets:

\begin{itemize}
\item the case $v\in V_{P_\varepsilon(T)}$ is not a concern since the cost doesn't 
change when considering $v\in V_T$ or $v\in V_{P_\varepsilon(T)}$; 

\item if $v \in U^T_{C_\varepsilon}$ then $v \in U^T_C$ or $D^T_C$; in the first case we ignore it, in the second case there is $\{v'\}=\max\{v''<v\mid v''\in V_{P_\varepsilon(T)}\}$, such that $v'\in D_C^T$. And so $cost_C(v)<cost_{C}(v')$;
\item if $v\in D^T_{C_\varepsilon}$, then $v\in D^T_{C}$ and $\#\Lambda_C(v)=0$; by construction $\varphi_{C_\varepsilon}(v)\leq \varphi_{C}(v)$. 
If $\varphi_{C_\varepsilon}(v)< \varphi_{C}(v)$ then also $\eta_{C_\varepsilon}$ is deleted and, since $f(\eta_{C_\varepsilon}(v))<f(v)$, $cost_C(v)<cost_C(\eta_{C_\varepsilon}(v))$.

So we are left with the case $\varphi_{C_\varepsilon}(v)= \varphi_{C}(v)$.
If $g(\eta_C(v))-f(v)>0.5\cdot \left(f(\varphi_{C}(v)-f(v))\right)$ then, again,
$cost_C(v)<cost_C(\eta_{C_\varepsilon}(v))$, for $f(\eta_{C_\varepsilon}(v))<f(v)$. 
Thus we always have $cost_C(v)\leq \max\{\varepsilon/2, cost_C(\eta_{C_\varepsilon}(v))\}$. 
\end{itemize}

Exchanging the role of $T$ and $G$ and repeating the same observations, we obtain that, if we consider $C\in\mathcal{C}^{o}(T,G)$, we have $d_I(T,G)\leq \parallel C \parallel_\infty\leq \max\{d_I(P_\varepsilon(T),P_\varepsilon(G)),\varepsilon/2\}$.

\hfill$ \blacksquare $

\bigskip\noindent
\underline{\textit{Proof of}  
\Cref{lem:extension_cost}.}

\smallskip\noindent

\begin{itemize}
\item if $v\leq r$ then $r'\geq \chi(v),\varphi(v)$ and $r'\geq \eta(v)$, so $cost(v)$ depends only on $C_{r,r'}$;
\item if $v\nleq r$ then it is either unused, if $v>r$, or it is deleted with $\#\Lambda(v)=0$. Then the cost of such deletion is either $0.5(f(v_r)-f(v))$ or $g(\eta(v))-f(v)$. If $C_{r,r'}\in\mathcal{C}^o_R(T_r,G_{r'})$, $g(\eta(v))=g_r$ and so we are done. Otherwise we simply have $g(\eta(v))\geq g_r$ and thus $g(\eta(v))-f(v)\geq g_r-f(v)$ implying the lower bound that we needed to prove.
\end{itemize}

\hfill$ \blacksquare $

\bigskip\noindent
\underline{\textit{Proof of}  
\Cref{prop:opt}.}

\smallskip\noindent

We just explore the different pieces of the cost function to assess the thesis:
\begin{itemize}
\item in the case of 
$(x,y)\in C(\mathcal{V})$, the cost of coupling $T_x$ and $G_y$ is given by $\sum_{y} a_{x,y}\parallel C_{x,y}\parallel _\infty$; 

\item if $u_x=1$ we have $\mid g(r_G)-f(x)\mid u_x$ which is the cost of deleting $x$ with $\#\Lambda(x)>1$. That is, any time $x\geq x'$ with $x'=\LCA(v,v')$ for $v,v'\in \pi_T(C(\mathcal{V}))$. Recall that, in this case, we have $\chi(x)\leq r_G$;

\item deleting $x$ with $\#\Lambda(x)=0$ is taken care by remaining part the cost function. 
For a vertex $x$ lets indicate with $x_f$ its father and 
 set $A_x=0.5(f(x_f)-f_x)(1-d_x)$ and $B_x = \lbrace \sum_y a_{v,y}\left( g_y-f_x \right) - K d_x\mid v<x_f \rbrace$. Now we try to unveil the meaning of $A_x$ and $B_x$.
Recall that deleting $x$ with $\#\Lambda(x)=0$ corresponds to having $d_x=0$. 
If $d_x =1$, then $A_x=0$ and $\max B_x<0$; while if $d_x >1$, both $A_x$ and $\max B_x$ are negative.  
Consider now $d_x=0$. 
Then $\varphi_{C_m}(x)= x'_f$, with $x'_f$ being father of some $x'\geq x$. Clearly $d_{x'}=0$ as well and $A_{x'}>A_x$.
In particular $A_{x'}=\max_{v\leq x'} 0.5 \left( f(\varphi_{C_m}(v))-f_v\right)$ (the same holds for $C^o$). 
Now, we turn to $B_x$. If $x<x'$ then $x_f<\varphi_{C^o}(x)=x'_f$, and so $d_{x_f}=0$, entailing $\max B_x=\min B_x=0$. Instead, if $x=x'$, we have by construction $v<x'_f$
with $a_{v,y}=1$. Then $\max \lbrace\sum_y a_{v,y} g_y\mid v<x_f \rbrace\geq g(\eta_{C^o}(x'))$ and so:
\begin{equation}
A_{x'}\leq \max_{v\leq x'}cost_{C^*}(v)\leq \max\lbrace A_{x'},\max B_{x'}\rbrace.
\end{equation}

\end{itemize}
\hfill$ \blacksquare $

\bigskip\noindent
\underline{\textit{Proof of}  
\Cref{cor:opt}.}

\smallskip\noindent

Given $C^*\in\mathcal{C}^*(T,G)$,
with $r=\LCA(\pi_T(C^*))$ and 
$r'=\LCA(\pi_G(C^*))$,
thanks to \Cref{cor:approx_1} we can approximate with $\Gamma^\uparrow$ and
$\Gamma^\downarrow$ the costs of the vertices in $T_r$ and $T_{r'}$. By \Cref{lem:extension_cost}, $H_{r,r'}$ then takes care of the vertices $v\nleq r$ and $w\nleq r'$.

\hfill$ \blacksquare $

\bibliography{references}

\end{document}